

\input epsf.tex

\def\2{{1\over 2}}

\def\d{\delta}
\def\a{\alpha}
\def\b{\beta}
\def\g{\gamma}

\def\s{\sigma}
\def\e{\epsilon}
\def\l{\lambda}
\def\o{\omega}

\def\fun#1#2#3{#1\colon #2\rightarrow #3}

\def\frac#1#2{{{#1} \over {#2}}}

\def\sqr{\sqrt}
\def\st{\;\colon\;}
\def\tends{\rightarrow}

\def\dr{ {\rm d} }

\def\C{{\bf C}}

\def\R{{\bf R}}
\def\N{{\bf N}}

\def\thm#1{\vskip 1 pc\noindent{\bf Theorem #1.\quad}\sl}
\def\lem#1{\vskip 1 pc\noindent{\bf Lemma #1.\quad}\sl}
\def\prop#1{\vskip 1 pc\noindent{\bf Proposition #1.\quad}\sl}

\def\proof{\rm\vskip 1 pc\noindent{\bf Proof.\quad}}
\def\fin{\par\hfill $\backslash\backslash\backslash$\vskip 1 pc}
\def\txt#1{\quad\hbox{#1}\quad}

\def\L{{\cal L}}

\def\s{\sigma}

\def\o{\omega}

\def\2{\frac{1}{2}}
\def\inn#1#2{{\langle #1 ,#2\rangle}}

\def\part{{\partial_{x}}}

\def\tr{{{}^{t}}}
\def\bc{{\cal B}}

\def\fc{{\cal F}}

\def\oc{{\cal O}}
\def\pc{{\cal P}}

\def\mc{{\cal M}}
\def\kc{{\cal K}}

\def\lk{{\Lambda^q}(\R^d)}



\baselineskip= 17.2pt plus 0.6pt
\font\titlefont=cmr17
\centerline{\titlefont Counting periodic orbits on fractals}
\vskip 1 pc
\centerline{\titlefont weighted by their Lyapounov exponents}
\vskip 4pc
\font\titlefont=cmr12
\centerline{         \titlefont {Ugo Bessi}\footnote*{{\rm 
Dipartimento di Matematica, Universit\`a\ Roma Tre, Largo S. 
Leonardo Murialdo, 00146 Roma, Italy.}}   }{}\footnote{}{
{{\tt email:} {\tt bessi@matrm3.mat.uniroma3.it} Work partially supported by the PRIN2009 grant "Critical Point Theory and Perturbative Methods for Nonlinear Differential Equations. Competing interests: the author declares none.}} 
\vskip 0.5 pc
 
\par
\vskip 2pc
\centerline{\bf Abstract}

Several authors have shown that Kusuoka's measure $\kappa$ on fractals is a scalar Gibbs measure; in particular, it maximises a pressure. There is also a different approach, in which one defines a matrix-valued Gibbs measure $\mu$ which induces both Kusuoka's measure $\kappa$ and Kusuoka's bilinear form. In the first part of the paper we show that one can define a "pressure" for matrix valued measures; this pressure is maximised by 
$\mu$. In the second part, we use the matrix-valued Gibbs measure $\mu$  to count periodic orbits on fractals, weighted by their Lyapounov exponents.

\vskip 2 pc
\centerline{\bf  Introduction}
\vskip 1 pc

We begin with a loose definition of the  fractals we consider. We are given $t$ affine, invertible contractions $\fun{\psi_i}{\R^d}{\R^d}$; the fractal $G$ is the unique compact set of $\R^d$ such that 
$$G=\bigcup_{i=1}^t\psi_i(G)  .  $$
As shown in [11], it is possible to define on $G$ a natural measure and bilinear form; these objects are connected ([6], [13], [16]) to Gibbs measures for systems of $d\times d$ matrices. We briefly outline the approach of [2] and [3], which derives Kusuoka's measure and bilinear form from a matrix-valued Gibbs measure.

First of all, under the hypotheses of Section 1 below the maps $\psi_i$ are the branches of the inverse of an expansive map $\fun{F}{G}{G}$, and the construction of Gibbs measures for expansive maps is a staple of dynamical systems theory ([12], [14], [23]). Let 
$q\in\{ 0,1,\dots,d \}$ and let $\Lambda^q(\R^d)$ be the space of $q$-forms on $\R^d$; since $\Lambda^q(\R^d)$ inherits an inner product from $\R^d$ we can define $M^q$, the space of self-adjoint operators (or symmetric matrices) on $\Lambda^q(\R^d)$. Let 
$$\fun{(D\psi_i)_\ast}{\Lambda^q(\R^d)}{\Lambda^q(\R^d)}$$
be the pull-back operator induced by the maps $D\psi_i$, i.e. 
$$[(D\psi_i)_\ast\o](v_1,\dots,v_q)=\o(D\psi_i\cdot v_1,\dots,D\psi_i\cdot v_k)  .  
\eqno (1)$$
The linear map $(D\psi_i)_\ast$ induces a push-forward operator
$$\fun{\Psi_i}{M^q}{M^q}  $$
by
$$(\Psi_i(A)\cdot v,w)=(A\cdot(D\psi_i)_\ast v,(D\psi_i)_\ast w)
\txt{or}\Psi_i(A)=\tr(D\psi_i)_\ast\cdot A\cdot(D\psi_i)_\ast  .  \eqno (2)$$
In the formula above, we have denoted by $(\cdot,\cdot)$ the inner product of 
$\Lambda^q(\R^d)$ and by $\tr B$ the transpose of the matrix $B$. 

For $\a\in(0,1]$ we define $C^{0,\a}(G,\R)$ in the standard way, i.e. $V\in C^{0,\a}(G,\R)$ if there is $C>0$ such that 
$$|V(x)-V(y)|\le C||x-y||^\a\txt{for all}x,y\in G  .  $$
For $V\in C^{0,\a}(G,\R)$ we define a Ruelle operator $\L_{G,V}$ in the following way: 
$$\fun{\L_{G,V}}{C(G,M^q)}{C(G,M^q)}$$
$$(\L_{G,V}A)(x)=\sum_{i=1}^t\Psi_i(A\circ{\psi_i(x)})\cdot e^{V\circ\psi_i(x)}  .  \eqno (3)$$
As we recall in Proposition 1.2 below, there is a unique $\b_{G,V}>0$ and a function 
$Q_{G,V}\in C(G,M^q)$ such that $Q_{G,V}(x)$ is positive-definite for all $x\in G$ and$$\L_{G,V}Q_{G,V}=\b_{G,V}Q_{G,V}  .  \eqno (4)$$
Moreover, the eigenvalue $\b_{G,V}$ is simple, i.e. $Q_{G,V}$ is unique up to multiplication by a scalar. 

By Riesz's representation theorem, the dual space of $C(G,M^q)$ is the space of $M^q$-valued Borel measures on $G$, which we denote by $\mc(G,M^q)$. Since $\L_{G,V}$ is a bounded operator by (3), its adjoint 
$$\fun{\L_{G,V}^\ast}{\mc(G,M^d)}{\mc(G,M^d)}  $$
is again bounded. Again by Proposition 1.2 below, there is $\mu_{G,V}\in\mc(G,M^d)$, unique up to multiplication by a scalar, such that for the same eigenvalue $\b_{G,V}$ of (4) we have 
$$\L_{G,V}^\ast\mu_{G,V}=\b_{G,V}\mu_{G,V}  .  \eqno (5)$$
The measure $\mu_{G,V}$ is semi-positive definite in the following sense: for all Borel sets 
$B\subset G$, $\mu_{G,V}(B)$ is a semi-positive definite matrix. 

Kusuoka's measure $\kappa_{G,V}$ is the scalar measure on $G$ defined by 
$$\kappa_{G,V}\colon=(Q_{G,V},\mu_{G,V})_{HS}  .  $$
On the right hand side of the formula above there is the Hilbert-Schmidt product  of the density $Q_{G,V}$ with the measure $\mu_{G,V}$; the details are in Section 1 below. An important fact is that $\kappa_{G,V}$ is ergodic for the expansive map $F$. 

In  dynamical systems, the positive eigenvector $\mu_{G,V}$ of $\L_{G,V}^\ast$ is usually called a Gibbs measure ([14], [23]) and that's how, in the following, we shall call our matrix-valued $\mu_{G,V}$; the logarithm of the eigenvalue $P(V)\colon=\log\b_{G,V}$  is called the pressure. This immediately raises the question whether our pressure, like its scalar counterpart, is the maximal value of a natural functional. The paper [5] gives a positive answer to this question, but [5] considers {\it scalar} Gibbs measures on a system of matrices, while in Theorem 1 below we are interested in defining a pressure for matrix-valued measures. 

When $M\in M^q$ is semi-positive definite, we write $M\ge 0$; we define $\kc_V$ as the set of the couples $\{ m,M \}$ where $m$ is a non-atomic, $F$-invariant probability measure on $G$ and $\fun{M}{G}{M^d}$ is a Borel function such that 
$$M(x)\ge 0
\txt{and}
(Q_{G,V}(x),M(x))_{HS}=1
\txt{for $m$-a.e.} x\in G  . \eqno (6)$$
We denote by $h(m)$ the entropy of $m$ with respect to $F$; the reader can find a definition in [24]. In Section 1 below we shall see that $\mu_{G,V}$  is absolutely continuous with respect to $\kappa_{G,V}$, and thus $\mu_{G,V}=M_{G,V}\cdot\kappa_V$, where $\fun{M_{G,V}}{G}{M^q}$ is a Borel function which satisfies (6), i.e. 
$\{ \kappa_{G,V},M_{G,V} \}\in\kc_V$.

\thm{1} Let the fractal $G$ satisfy hypotheses (F1)-(F4) and $(ND_q)$ of Section 1 below. Let $V\in C^{0,\a}(G,\R)$ and let $\kappa_{G,V}=(Q_{G,V},\mu_{G,V})_{HS}$ be Kusuoka's measure for the potential $V$. Let $\b_{G,V}$ be the eigenvalue of $\L_{G,V}$ as illustrated in (4) and let us set 
$P(V)=\log\b_{G,V}$. Then, 
$$P(V)=h(\kappa_{G,V})+\int_G V(x)\dr\kappa_{G,V}(x)+$$
$$\sum_{i=1}^t\int_G\log(Q_{G,V}(x),\tr\Psi_i\circ M_{G,V}\circ F(x))_{HS}
\dr\kappa_{G,V}(x) = $$
$$\sup_{\{ m,M \}\in\kc_V}\Bigg\{
h(m)+\int_G V(x)\dr m(x)+$$
$$\sum_{i=1}^t\int_{G}\log(Q_{G,V} (x),\tr\Psi_i\circ M\circ F (x))_{HS}
\dr m(x)
\Bigg\}   .  \eqno (7)  $$

\rm

\vskip 1pc

Some remarks. First, in the classical case of [15] and [24] the topological entropy (i.e. the pressure of the function $V\equiv 0$) is always positive, while in our case it can be positive, negative or zero. Indeed, the contribution of the last term in (7) is non-positive, since the maps $\psi_i$ are contractions. It is well known (see for instance [7]) that for a famous example, the harmonic Sierpinski gasket with $q=1$, the pressure of the zero function is negative, since $\b_{G,0}=\frac{3}{5}$. 

The second remark is that the raison d'\^etre of hypotheses (F1)-(F4) is to translate Theorem 1 to a theorem on the shift on $t$ symbols, where we shall prove it. We have stated it on the fractal to underline the existence of a Riemannian structure on $G$ which is natural with respect to the dynamics. 

The last remark is that the measure $\kappa_{G,V}$ coincides with that found in [5] and [13];  this is shown in the appendix of [13] and we give another proof after Lemma 1.3 below. We also remark that the papers [5], [13] and [16] are more general than ours: they consider the pressure of a real power $s$ of the Lyapounov exponent, while we fix $s=2$. 

Another property of the scalar Gibbs measures is that they satisfy a central limit theorem; we recall its statement from [23] and [14]. 

Let $m$ be a $F$-invariant probability measure on $G$ and let $\phi$ be any function in 
$C^{0,\a}(G,\R)$ which satisfies 
$$\int_G\phi(x)\dr m(x)=0   .     \eqno (8)$$
We define 
$$\s^2=\int_G\phi^2(x)\dr m(x)+
2\sum_{j\ge 1}\int_G\phi(x)\cdot\phi\circ F^j(x)\dr m(x)  .  $$
We say that $m$ satisfies a central limit theorem if, for all functions $\phi$ as above, the  three points below hold. 

\noindent 1) $\s\in[0,+\infty)$. 

\noindent 2) If $\s=0$, then $\phi= u\circ F-u$ for some 
$u\in L^2(G,m)$. 

\noindent 3) If $\s>0$ then for all intervals $A\subset\R$ we have that 
$$m\left\{
x\in G\st\frac{1}{\sqrt n}\sum_{j=0}^{n-1}\phi\circ F^j(x)\in A
\right\}   \tends\frac{1}{\s\sqrt{2\pi}}\int_Ae^{-\frac{z^2}{2\s^2}}\dr z  .  $$

\vskip 1pc

It turns out ([23]) that $m$ satisfies a central limit theorem if it has an exponential decay of correlations, which means the following: there are $C,\d>0$ such that, for all 
$\phi,\psi\in C^{0,\a}(G,\R)$, we have 
$$\left\vert
\int_G\phi(x)\psi\circ F^n(x)\dr m(x)-
\int_G\phi(x)\dr m(x)\int_G\psi(x)\dr m(x)
\right\vert  \le C e^{-n\d}  .  $$
Moreover, the constant $C$ only depends on $||\psi||_{L^1}$ and 
$||\phi||_{C^{0,\a}}+||\phi||_{L^1}$, while $\d$ doesn't depend on anything.

Since Kusuoka's measure $\kappa_{G,V}$ satisfies the formula above by point 6) of Proposition 1.2 below, we have that $\kappa_{G,V}$ satisfies a central limit theorem. 

Next, we want to count periodic orbits, a task which does not look very interesting because our hypotheses on the fractal $G$ and the expansive map $F$ imply almost immediately that there is a bijection between the periodic orbits of $F$ and those of a Bernoulli shift on $t$ elements. What we are going to do is to count periodic orbits, but weighted by their Lyapounov exponents. In order to be more precise, we need some notation. 

Let 
$$x_0=x,\quad x_1=F(x_0),\quad\dots,\quad x_{i+1}=F(x_i),\quad\dots  $$
be a periodic orbit with minimal period $n$; in other words, $x_i=x_{i+n}$ for all $i\ge 0$ and $n$ is the smallest integer with this property. We group in a set $\tau$ all the points of the orbit 
$$\tau=\{ x_0,x_1,\dots,x_{n-1} \}  \eqno (9)$$
and we define 
$$\l(\tau)=n=\sharp\tau,\qquad
\Psi_{x_0\dots x_{n-1}}=\Psi_{x_{n-1}}\circ\dots\circ\Psi_{x_0},\qquad
\tr\Psi_{x,n}=\tr\Psi_{x_0\dots x_{n-1}}=\tr\Psi_{x_0}\circ\dots\circ\tr\Psi_{x_{n-1}}  .  \eqno (10)$$
Note the reverse order in the definition of $\Psi_{x_0\dots x_{n-1}}$.

Let $\hat V\in C^{0,\a}(G,\R)$ be positive and let us set 
$\hat V(\tau)=V(x_0)+\dots +V(x_{n-1})$; for $r\ge 1$ we define, for $c\in\R$ such that 
$P(-c\hat V)=0$,
$$\pi(r)=\sum_{\tau\st e^{|c|\hat V(\tau)}\le r}
{\rm tr}[\tr\Psi_{x_0}\circ\dots\circ\tr\Psi_{x_{n-1}}]  
\eqno (11)$$
where $\tau$ and $x_0\dots x_{n-1}$ are related as in (9); by the cyclical invariance of the trace of a product, ${\rm tr}[\tr\Psi_{x_0}\circ\dots\circ\tr\Psi_{x_{n-1}}]$ does not depend on where we place $x_0$ along the periodic orbit. As we shall see in Section 3, the summands on the right side of (11) are all positive. 

We saw in Theorem 1 that the topological entropy (i.e., $P(V)$ with $V\equiv 0$) can be positive, negative or zero; we have to exclude the case in which it is zero. Moreover, we have to make the usual distinction ([14], [15]) between the weakly mixing and non-weakly mixing case, because they have different asymptotics. For simplicity, in point 1) of Theorem 2 below we shall consider only one case of a non-weakly mixing suspension, i.e. the one with $\hat V\equiv 1$. In the scalar case, a delicate study of the dynamics (chapters 4-6 of [14]) shows the topological weak mixing property for the suspension by the function $\hat V$ is equivalent to the fact that the zeta function can be extended to a non-zero analytic function on ${\rm Re}(s)=1$, save for a simple pole at $1$. In point 2) of Theorem 2 below we shall simply suppose that this extension is possible, since we don't know how to adapt [14] to our situation. 

\thm{2} Let $\hat V\in C^{0,\a}(G,\R)$ be a potential such that $\hat V(x)>0$ for all 
$x\in G$, and let $c$ be the unique real such that the pressure of $-c\hat V$ is zero. Let 
$\pi$ be defined as in (11) and let us suppose that $c\not=0$. 

\noindent 1) Let $\hat V\equiv 1$ and let us set for simplicity $\b=\b_{G,0}$. Then, if we  take $c=\log\b$ in (11) we have 
$$\limsup_{r\tends+\infty}\frac{\pi(r)\log r}{r}\le\frac{\b}{\b-1}\log\b     \eqno (12)$$
if $\b>1$; if $\b<1$, $\pi(r)$ is bounded. 

\noindent 2) Let now $\hat V\in C^{0,\a}(G,\R)$ be strictly positive and let us suppose that the zeta function of Section 4 below can be extended to a continuous, non-zero function on $\{ {\rm Re}(s)\ge 1 \}\setminus \{ 1 \}$, with a simple pole in $s=1$ in the case $c>0$. Then, 
$$\limsup_{r\tends+\infty}\frac{\pi(r)\log r}{r}\le
\left\{
\eqalign{
1&\txt{if}c>0\cr
0&\txt{if}c<0   . 
}
\right.        \eqno (13)$$
Naturally, in the second case the inequality is an equality.

\rm

\vskip 1pc

Note that in the case $c>0$ we are only giving an estimate from above on the asymptotics, while [14] gives much more, an equality. 

As a last remark, we make no pretence at novelty: already in [22] D. Ruelle considered transfer operators and zeta functions acting on matrices. Moreover, the paper [10] applies similar techniques to a different problem, the Weyl statistics for the eigenvalues of the Laplacian on fractals. 

This paper is organised as follows. In Section 1 we set the notation and define our family of fractals; it is the standard definition of an iterated function system. We end this section recalling the construction of Kusuoka's measure $\kappa_{\hat V}$ from [2] and [3]. In Section 2 we prove Theorem 1, in Section 3 we recall an argument of [14] for the proof of Theorem 2 in the non-weakly mixing case; we shall apply this argument to the simple case $\hat V\equiv 0$ (or $\hat V\equiv 1$, since all constants induce the same Kusuoka's measure.) In Section 4 we prove Theorem 2 in the "weakly mixing" case; we follow closely the arguments of [14]. 

\vskip 2pc

\centerline{\bf \S 1}
\centerline{\bf Preliminaries and notation}
\vskip 1pc

\noindent{\bf The Perron-Frobenius theorem.} We recall without proof a few facts from [23]; we refer the reader to [4] for the original treatment. 

Let $E$ be a real vector space; a cone on $E$ is a subset $C\subset E\setminus \{ 0 \}$ such that, if $v\in C$, then $tv\in C$ for all $t>0$. We define the closure $\bar C$ of 
$C$ as the set of all $v\in E$ such that there is $w\in C$ and $t_n\searrow 0$ such that 
$v+t_nw\in C$ for all $n\ge 1$. 

In the following, we shall always deal with convex cones $C$ such that 
$$\bar C\cap(-\bar C)=\{ 0 \}  .  \eqno (1.1)$$
Given $v_1,v_2\in C$ we define
$$\a(v_1,v_2)=\sup\{
t>0\st v_2-tv_1\in C
\}  ,    $$
$$\frac{1}{\b(v_1,v_2)}=\sup\{
t>0\st v_1-tv_2\in C
\}   $$
and 
$$\theta(v_1,v_2)=
\log\frac{\b(v_1,v_2)}{\a(v_1,v_2)}  . $$
We define an equivalence relation on $C$ saying that $v\simeq w$ if $v=tw$ for some 
$t>0$. It turns out that $\theta$ is a metric on $\frac{C}{\simeq}$, though it can assume the value $+\infty$. As shown in [23], it separates points thanks to (1.1). 

We recall the statement of the Perron-Frobenius theorem. 

\thm{1.1} Let $C\subset E$ be a convex cone satisfying (1.1) and let $\fun{L}{E}{E}$ be a linear map such that $L(C)\subset C$. Let us suppose that 
$$D\colon=\sup\{
\theta(L(v_1),L(v_2))\st v_1,v_2\in C
\}  <+\infty  .  $$
Then, the following three points hold. 

\noindent 1) For all $v_1,v_2\in C$ we have that 
$$\theta(L(v_1),L(v_2))\le(1-e^{-D})\theta(v_1,v_2)  .  $$

\noindent 2) If $(\frac{C}{\simeq},\theta)$ is complete, then $L$ has a unique fixed point 
$v\in\frac{C}{\simeq}$. Since $L$ is linear, this means that there is $v\in C$, unique up to multiplication by a positive constant, and a unique $\l>0$ such that 
$$Lv=\l v  .  $$

\noindent 3) Let $(\frac{C}{\simeq},\theta)$ be complete and let $v$ be the fixed point of  point 2). Then, if $v_0\in C$ we have that 
$$\theta(L^n v_0,v)\le
\frac{(1-e^{-D})^n}{e^{-D}}\theta(v_0,v)  .  $$

\rm

\vskip 1pc

\noindent{\bf Fractal sets.} The fractals we consider are a particular case of what in [9] (see definitions 1.3.4 and 1.3.13) are called "post-critically finite self-similar structures."

\noindent {\bf(F1)} We are given $t$ affine bijections $\fun{\psi_i}{\R^d}{\R^d}$, 
$i\in\{ 1,\dots,t \}$, satisfying 
$$\eta\colon=\sup_{i\in\{ 1,\dots,t \}}Lip(\psi_i)<1  .   \eqno (1.2)$$
By Theorem 1.1.7 of [9], there is a unique non empty compact set $G\subset\R^d$ such that
$$G=\bigcup_{i\in\{ 1,\dots,t \}}\psi_i(G) . \eqno (1.3)$$
In the following, we shall always rescale the norm of $\R^d$ in such a way that 
$${\rm diam}(G)\le 1   .  \eqno (1.4)$$

If (F1) holds, then the dynamics of $F$ on $G$ can be coded. Indeed, we define 
$\Sigma$ as the space of sequences 
$$\Sigma=
\{  
( x_i )_{i\ge 0}\st x_i\in(1,\dots,n) 
\}    $$
with the product topology. One metric which induces this topology is the following: for 
$\g\in(0,1)$, we set 
$$d_\g(( x_i )_{i\ge 0},( y_i )_{i\ge 0})=
\inf\{
\g^{k+1}\st x_i=y_i\txt{for}i\in\{ 0,\dots,k \}
\}    \eqno (1.5)$$
with the convention that the $\inf$ on the empty set is $1$.

We define the shift $\s$ as 
$$\fun{\s}{\Sigma}{\Sigma},\qquad
\fun{\sigma}{( x_0,x_1,x_2,\dots )}{( x_1,x_2,x_3,\dots )}   .   $$
If $x=(x_0x_1\dots)\in\Sigma$, we set $ix=(ix_0x_1\dots)$.  

For $x=(x_0x_1\dots)\in\Sigma$ we define the cylinder of $\Sigma$ 
$$[x_0\dots x_l]=\{
( y_i )_{i\ge 0}\in\Sigma\st y_i=x_i\txt{for}i\in(0,\dots,l)
\}  .  $$
For $x=(x_0x_1\dots)\in\Sigma$ and $l\in\N$ we set 
$$\psi_{x_0\dots x_l}=\psi_{x_0}\circ\dots\circ\psi_{x_l}  $$
and
$$[x_0\dots x_l]_G=\psi_{x_0}\circ\psi_{x_{1}}\circ\dots\circ\psi_{x_l}(G)  .  \eqno (1.6)$$
Formula (1.6) immediately implies that, for all $i\in(1,\dots,t)$, 
$$\psi_i([x_1\dots x_l]_G)=[i x_1\dots x_l]_G  .  \eqno (1.7)$$
Since the maps $\psi_i$ are continuous and $G$ is compact, the sets 
$[x_0\dots x_l]_G\subset G$ are compact. By (1.3) we have that 
$\psi_i(G)\subset G$ for $i\in(1,\dots,t)$; together with (1.6) this implies that, for all 
$( x_i )_{i\ge 0}\in\Sigma$ and $l\ge 1$,
$$[x_0\dots x_{l-1} x_l]_G\subset [x_0\dots x_{l-1}]_G  .  \eqno (1.8)$$
From (1.2), (1.4) and (1.6) we get that, for all $l\ge 0$, 
$${\rm diam} ([x_0\dots x_l]_G)\le\eta^{l+1}   .   \eqno (1.9)$$
By (1.8), (1.9) and the finite intersection property we get that, if 
$( x_i )_{i\ge 0}\subset\Sigma$, then 
$$\bigcap_{l\ge 1}[x_0\dots x_{l}]_G$$
is a single point, which we call $\Phi(( x_i )_{i\ge 0})$; Formula (1.9) implies in a standard way that the map $\fun{\Phi}{\Sigma}{G}$ is continuous. It is not hard to prove, using (1.3), that $\Phi$ is surjective. In our choice of the metric $d_\g$ on $\Sigma$ we shall always take $\g\in(\eta,1)$; we endow $G$ with the Euclidean distance on $\R^d$. With this choice of the metrics, (1.5) and (1.9) easily imply that $\Phi$ is 1-Lipschitz. By (1.7) and the definition of $\Phi$ we see that  
$$\psi_i\circ\Phi(x)=\Phi(ix)   .   \eqno (1.11)$$

\noindent {\bf (F2)} We ask that, if $i\not=j$, then $\psi_i(G)\cap\psi_j(G)$ is a finite set. We set 
$$\fc\colon=\bigcup_{i\not =j}\psi_i(G)\cap\psi_j(G)  .  $$

\noindent {\bf (F3)} We ask that, if $i=(i_0i_1\dots)\in\Sigma$ and $l\ge 1$, then 
$\psi_{i_0\dots i_{l-1}}^{-1}(\fc)\cap\fc=\emptyset$.

We briefly prove that (F1)-(F3) imply that the coding $\Phi$ is finite-to-one and that the set 
$N\subset\Sigma$  where $\Phi$ is not injective is countable. More precisely, we assert that the points of $G$ in  
$$H\colon=\fc\cup\bigcup_{l\ge 0}\bigcup_{x\in\Sigma}
\psi_{x_0\dots x_l}(\fc)  \eqno (1.12)$$
have at most $t$ preimages, and that those outside $H$ have only one preimage. 

We begin by showing this latter fact. If $y\in G\setminus\fc$, then $y=\Phi(( x_j )_{j\ge 0})$ can belong by (F2) to at most one $\psi_i(G)$, and thus there is only one choice for $x_0$. Using the fact that $\psi_{x_0}$ is a diffeomorphism, we see that there is at most one choice for $x_1$ if moreover $y\not\in\psi_{x_0}(\fc)$; iterating, we see that the points in $H^c$ have one preimage. 

If $y\in\fc$, then $y$ can belong to at most $t$ sets $\psi_i(G)$, and thus there are at most $t$ choices for $x_0$. As for the choice for $x_1$, once we have fixed $x_0$ then we have that $y\in\psi_{x_0}\circ\psi_{x_1}(G)$, i.e. that $\psi_{x_0}^{-1}(y)\in\psi_{x_1}(G)$. Since $y\in\fc$, (F3) implies that $\psi_{x_0}^{-1}(y)\not\in\fc$ and thus we have at most one choice for $x_1$. Iterating, we see that there is at most one choice for $x_2,x_3$, etc... If 
$y\not\in\fc$ but $y\in\psi_{x_0}(\fc)$, then we have one choice for $x_0$ but $t$ choices for $x_1$; arguing as above we see that there is at most one choice for $x_2$, 
$x_3$, etc.... Iterating, we get the assertion. 

\noindent {\bf (F4)} We ask that there are disjoint open sets $\oc_1,\dots,\oc_n\subset\R^d$ such that
$$G\cap\oc_i=\psi_i(G)\setminus\fc\txt{for}i\in(1,\dots,n)  .  $$
We define a map $\fun{F}{\bigcup_{i=1}^n\oc_i}{\R^d}$ by
$$F(x)=\psi_i^{-1}(x)  \txt{if}x\in\oc_i  .  $$
This implies the second equality below; the first one holds by the formula above if we ask that $\psi(\oc_i)\subset\oc_i$.    
$$F\circ\psi_i(x)=x\qquad\forall x\in\oc_i 
\txt{and}
\psi_i\circ F(x)=x\qquad\forall x\in\oc_i  .  \eqno (1.13)$$
We call $a_{i}$ the unique fixed point of $\psi_i$; note that, by (1.3), $a_{i}\in G$. If 
$x\in\fc$, we define $F(x)=a_{i}$ for some arbitrary $a_{i}$.  This defines $F$ as a Borel map on all of $G$, which satisfies (1.13).

From (1.6) and the definition of $\Phi$ we see that $\Phi([x_0\dots x_l])=[x_0\dots x_l]_G$. We saw above that $\Phi$ is finite-to-one and that there are at most countably many points with more that one preimage; an immediate consequence is that the set 
$$\Phi^{-1} ([x_0\dots x_l]_G)\setminus [x_0\dots x_l]     $$
is at most countable. Let us suppose that $\Phi(x_0x_1\dots)\in\oc_i$; the discussion after (F3) implies that there is only one choice for $x_0$ and it is $x_0=i$; by the first formula of (1.13) this implies the second equality below; the first and last ones are the definition of 
$\Phi$. 
$$F\circ\Phi(x_0x_1\dots)=F\left(
\bigcap_{l\ge 0}\psi_{x_0x_1\dots x_l}(G)
\right) =
\bigcap_{l\ge 0}\psi_{x_1x_2\dots x_l}(G)=\Phi\circ\s(x_0x_1\dots)   .   $$
This yields the first equality below, while the second one is (1.11). 
$$\left\{\matrix{
\Phi\circ\s(x)=F\circ\Phi(x)\txt{save possibly when $\Phi(x)\in \fc$,}\cr
\Phi(ix) =\psi_i(\Phi(x))   \qquad\forall x\in\Sigma,
\quad\forall i\in (1,\dots,n)   .   
}
\right.       \eqno (1.14)  $$
In other words, up to a change of coordinates, shifting the coding one place to the left is the same as applying $F$.

A particular case we have in mind is the harmonic Sierpinski gasket on $\R^2$ (see for instance [7]). We set
$$T_1=\left(
\matrix{
\frac{3}{5},&0\cr
0,&\frac{1}{5}
}
\right)  ,  \quad
T_2=\left(
\matrix{
\frac{3}{10},&\frac{\sqr 3}{10}\cr
\frac{\sqr 3}{10},&\frac{1}{2}
}
\right)  ,\quad
T_3=\left(
\matrix{
\frac{3}{10},&-\frac{\sqr 3}{10}\cr
-\frac{\sqr 3}{10},&\frac{1}{2}
}
\right)   ,  $$
$$A=\left(
\matrix{
0\cr
0
}
\right),\qquad
B=\left(
\matrix{
1\cr
\frac{1}{\sqrt 3}
}
\right) ,\qquad
C=\left(
\matrix{
1\cr
-\frac{1}{\sqrt 3}
}
\right)  $$
and 
$$\psi_1(x)=T_1(x),\quad \psi_2(x)=
B   +T_2\left(
x-B
\right),\quad
\psi_3(x)=
C   +T_3\left(
x-C
\right)  .  $$
Referring to Figure 1 below, $\psi_1$ brings the triangle $G_0\colon=ABC$ into $Abc$; 
$\psi_2$ brings $G_0$ into $Bac$ and $\psi_3$ brings $G_0$ into $Cba$. We set 
$\fc=\{ a,b,c \}$ and take $\oc_1$, $\oc_2$, $\oc_3$ as three disjoint open sets which contain, respectively, the triangle $Abc$ minus $b,c$, $Bca$ minus $c,a$ and $Cba$ minus $a,b$; in this way hypotheses (F1)-(F4) hold. It is easy to check that also hypothesis $(ND)_1$ below is satisfied. 

We define the map $F$ as 
$$F(x)=\psi_i^{-1}(x)
\txt{if}
x\in\oc_i $$
and we extend it arbitrarily on $\fc=\{ a,b,c \}$, say $F(a)=A$, $F(b)=B$ and $F(c)=C$. 



For $q\in(1,\dots,d)$ we consider $\lk$, the space of $q$-forms on $\R^d$; if $(\cdot,\cdot)$ is an inner product on $\R^d$, it induces an inner product on the monomials of 
$\Lambda^q(\R^d)$ by 
$$((v_1\wedge\dots\wedge v_q),(w_1\wedge\dots\wedge w_q))=
(v_1\wedge\dots\wedge v_q) (i(w_1),\dots, i(w_q))       $$
where $\fun{i}{\Lambda^1(\R^d)}{\R^d}$ is the Riesz map, i.e. the natural identification of a Hilbert space with its dual. This product extends by linearity to all $\Lambda^q(\R^d)$ which thus becomes a Hilbert space. 

As a consequence, if $A$ is a matrix bringing $\Lambda^q(\R^d)$ into itself, we can define its adjoint $\tr A$. We shall denote by $M^q$ the space of all self-adjoint operators (a.k.a. symmetric matrices) on $\lk$. Again since $\lk$ is a Hilbert space, we can define the trace ${\rm tr}(A)$ of $A\in M^q$ and thus the Hilbert-Schmidt product 
$$\fun{(\cdot,\cdot)_{HS}}{M^q\times M^q}{\R}$$
$$(A,B)_{HS}\colon={\rm tr}(A\tr B)={\rm tr}(AB) $$
where the second equality comes from the fact that $B$ is symmetric. 

Thus, $M^q$ is a Hilbert space whose norm is  
$$||A||_{HS}=\sqrt{(A,A)_{HS}}  .  $$
Since $\Lambda^k(\R^d)$ is a Hilbert space, we shall see its elements $b$ as column vectors, on which the elements $M\in M^q$ act to the left: $Mb$. 

\noindent{\bf Matrix-valued measures. } Let $(E,\hat d)$ be a compact metric space; in the following, $(E,\hat d)$ will be either one of the spaces $(\Sigma,d_\g)$ or $(G,||\cdot||)$ defined above. We denote by $C(E,M^q)$ the space of continuous functions from $E$ to $M^q$ and by $\mc(E,M^q)$ the space of the Borel measures on $E$ valued in $M^q$. If 
$\mu\in\mc(E,M^q)$, we denote by $||\mu||$ its total variation measure; it is standard ([20]) that $||\mu||(E)<+\infty$ and that 
$$\mu=M\cdot||\mu||$$
where $\fun{M}{E}{M^q}$ is Borel and 
$$||M(x)||_{HS}=1\txt{for $||\mu||$-a.e.} x\in E  .  \eqno (1.15)$$
Let now $\mu\in\mc(E,M^q)$ and let $\fun{f}{E}{M^q}$ be a Borel function such that 
$||f||_{HS}\in L^1(E,||\mu||)$. We define 
$$\int_E(f,\dr\mu)_{HS}\colon=
\int_E(f(x),M(x))_{HS}\dr||\mu||(x)  .  \eqno (1.16)$$
Note that the integral on the right converges: indeed, the function $x\rightarrow(f(x),M(x))$ is in $L^1(G,||\mu||)$ by our hypotheses on $f$, (1.15) and Cauchy-Schwarz. A particular case we shall use often is when $f\in C(E,M^q)$. 

If $\mu\in\mc(E,M^q)$ and $\fun{a,b}{E}{\lk}$ are Borel functions such that 
$||a||\cdot||b||\in L^1(E,||\mu||)$, we can define 
$$\int_E(a,\dr\mu\cdot b)\colon=
\int_E(
a(x),M(x)b(x)
)\dr||\mu||(x)  .  $$
Note again that the function $x\rightarrow(a(x),M(x)b(x))$ belongs to $L^1(G,||\mu||)$ by our hypotheses on $a$, $b$, (1.15) and Cauchy-Schwarz. We have written $M(x)b(x)$ because, as we said above, we consider the elements of $\Lambda^q(\R^d)$ as column vectors. Since 
$$(a\otimes b,M)_{HS}=(a,Mb)  $$
we have that 
$$\int_E(a,\dr\mu\cdot b)=
\int_E(a\otimes b,\dr\mu)_{HS}  .    \eqno (1.17)$$

By Riesz's representation theorem, the duality coupling $\inn{\cdot}{\cdot}$ between $C(E,M^q)$ and $\mc(E,M^q)$ is given by 
$$\inn{f}{\mu}\colon=
\int_E(f,\dr\mu)_{HS}   $$
where the integral on the right is defined by (1.16). 

Let $\mu\in\mc(E,M^q)$; we say that $\mu\in\mc_+(E,M^q)$ if $\mu(A)$ is a semi-positive-definite matrix for all Borel sets $A\subset E$. By Lusin's theorem (see [2] for more details), this is equivalent to saying that 
$$\inn{f}{\mu}\ge 0  \eqno (1.18)$$
for all $f\in C(E,M^q)$ such that $f(x)\ge 0$ (i.e. it is positive-semidefinite) for all $x\in E$. An equivalent characterisation is that, in the polar decomposition (1.15), $M(x)\ge 0$ for 
$||\mu||$-a.e. $x\in E$. 

As a consequence of the characterisation (1.18), $\mc_+(E,M^q)$ is a convex set closed for the weak$\ast$ topology of $\mc(E,M^q)$. 

If $\mu\in\mc(E,M^q)$ and $\fun{A}{E}{M^q}$ belongs to $L^1(E,||\mu||)$, we define the scalar measure $(A,\mu)_{HS}$ by 
$$\int_E f(x)\dr(A,\mu)_{HS}(x)=
\int_E(f(x)A(x),\dr\mu(x))_{HS}$$
for all $f\in C(E,\R)$. Recalling the definition (1.16) of the integral on the right we have that 
$$(A,\mu)_{HS}=(A,M)_{HS}\cdot||\mu||  $$
where $\mu=M\cdot||\mu||$ is the polar decomposition of $\mu$ as in (1.15). 

We recall (a proof is, for instance, in [2]) that, if $A,B,C\in M^q$, if $A\ge B$ and
$C\ge 0$, then 
$$(A,C)_{HS}\ge(B,C)_{HS}   .    \eqno (1.19)$$
As a consequence we have that, if $f,g\in C(E,M^q)$ satisfy $f(x)\ge g(x)$ for all $x\in E$  and if $\mu\in\mc_+(E,M^q)$, then 
$$\int_E(f(x),\dr\mu(x))_{HS}\ge \int_E(g(x),\dr\mu(x))_{HS}  .  \eqno (1.20)$$

Let $Q\in C(E,M^q)$ be such that $Q(x)>0$ for all $x\in E$; by compactness we can find 
$\e>0$ such that
$$\e Id\le Q(x)\le\frac{1}{\e}Id\qquad\forall x\in E  .  \eqno (1.21)$$
It is shown in [2] that there is $C_1=C_1(\e)>0$ such that, if $Q$ satisfies (1.21) and 
$\mu\in\mc_+(E,M^q)$, then 
$$\frac{1}{C_1(\e)}||\mu||\le
(Q,\mu)_{HS}\le
C_1(\e)||\mu||  .  \eqno (1.22)$$
Let $Q\in C(E,M^q)$ satisfy (1.21); we say that $\mu\in\pc_Q$ if $\mu\in\mc_+(E,M^q)$ and 
$$\inn{Q}{\mu}=1  .  \eqno (1.23)$$
We saw after Formula (1.18) that $\mc_+(E,M^q)$ is convex and closed for the weak$\ast$ topology; Formula (1.23) implies that $\pc_Q$ is also closed and convex; (1.22) and (1.23) imply that it is compact. 

\noindent{\bf Cones in function spaces.} We say that $f\in C_+(E,M^q)$ if $f\in C(E,M^q)$ and $f(x)$ is positive-definite for all $x\in E$. It is easy to see that $C_+(E,M^q)$ is a convex cone in $C(E,M^q)$ satisfying (1.1). We let $\theta^+$ denote the hyperbolic distance on $\frac{C_+(E,M^q)}{\simeq}$. 

Let $a>0$ and $\a\in(0,1]$; we say that $A\in LogC^{a,\a}_+(E)$ if $A\in C_+(E,M^q)$ and for all $x,y\in E$ we have 
$$A(y)e^{-a\hat d(x,y)^\a}\le A(x)\le A(y)e^{a\hat d(x,y)^\a}   $$
where the inequality is the standard one between symmetric matrices. 

It is easy to check that $LogC^{a,\a}_+(E)$ is a convex cone satisfying (1.1). We let 
$\theta^a$ denote the hyperbolic distance on $\frac{LogC^{a,\a}_+(E)}{\simeq}$.

Next, we introduce the subcone of the elements $A\in LogC^{a,\a}_+(E)$ whose eigenvalues are bounded away from zero. In other words, we define 
$$||A||_{\infty}=\sup_{x\in E}||A(x)||_{HS}$$
and for $\e>0$ we set 
$$LogC^{a,\a}_\e(E)=\{
A\in LogC_+^{a,\a}(E)\st A(x)\ge\e||A||_\infty Id\qquad\forall x\in E
\}   .  $$

\noindent{\bf The Ruelle operator.} Recall that the maps $\psi_i$ of Section 1 are affine and injective; in particular, $D\psi_i$ is a constant, nondegenerate matrix. We define the pull-back $\fun{(D\psi_i)_\ast}{\Lambda^q(\R^d)}{\Lambda^q(\R^d)}$ as in (1) and the push-forward $\fun{\Psi_i}{M^q}{M^q}$ as in (2). Let $\a\in(0,1]$ and let $V\in C^{0,\a}(G,\R)$.  

The Ruelle operator $\L_{G,V}$ on $C(G,M^q)$ is the one defined by (3).

We also define a Ruelle operator $\L_{\Sigma,V}$ on $C(\Sigma,M^q)$: defining $(ix)$ as in Section 1 we set 
$$\fun{\L_{\Sigma,V}}{C(\Sigma,M^q)}{C(\Sigma,M^q)}$$
$$(\L_{\Sigma,V} A)(x)=\sum_{i=1}^t e^{V\circ\Phi(ix)}\Psi_i(A(ix))  .  
\eqno (1.24)$$
Note that $V\circ\Phi$ is $\a$-H\"older: indeed, $V$ is $\a$-H\"older and we saw in Section 1 that $\Phi$ is 1-Lipschitz.

In order to apply the Perron-Frobenius theorem, we need a nondegeneracy hypothesis. 

\vskip 1pc

\noindent{\bf (ND)${}_q$} We suppose that there is $\g>0$ such that, if $c,e\in\lk$, we can find 
$i\in(1,\dots,t)$ such that 
$$|((D\psi_{i})_\ast c, e)|\ge\g ||c||\cdot||e||  .  $$
In other words, we are asking that, for all $c\in\Lambda^q(\R^d)\setminus\{ 0 \}$, 
$(D\psi_{i})_\ast c$ is a base of $\Lambda^q(\R^d)$; clearly, this implies that $t\ge\left(\matrix{d\cr q}\right)$. 

\vskip 1pc

We state the standard theorem on the existence of Gibbs measures, whose proof ([14], [23]) adapts easily to our situation; the details are in [2] and [3]. 

\prop{1.2} Let $(E,\hat d,S)$ be either one of $(\Sigma,d_\g,\s)$ or $(G,||\cdot||,F)$; let 
$V\in C^{0,\a}(\Sigma,\R)$ such that $V=\tilde V\circ\Phi$ for some  
$\tilde V\in C^{0,\a}(G,\R)$. Let the Ruelle operators $\L_{\Sigma,V}$ and $\L_{G,V}$ be defined as in (1.24) and (3) respectively. Then, the following holds. 

\noindent 1) There are $a,\e>0$, $Q_{E,V}\in LogC^{a-1,\a}_{\e}(E)$ and $\b_{E,V}>0$ such that 
$$\L_{E,V} Q_{E,V}=\b_{E,V} Q_{E,V}  .  $$
If $Q^\prime \in C(E,M^q)$ is such that 
$$\L_{E,V} Q^\prime=\b_{E,V} Q^\prime   ,   $$
then $Q^\prime=\eta Q_{E,V}$ for some $\eta\in\R$.  

\noindent 2) Since $Q_{E,V}\in LogC^{a-1,\a}_{\e(a)}(E)$, $Q_{E,V}$ satisfies (1.21) and we can define $\pc_{Q_{E,V}}$ as in (1.23). Then, there is a unique 
$\mu_{E,V}\in\pc_{Q_{E,V}}$ such that 
$$\L_{E,V}^\ast\mu_{E,V}=\b_{E,V}\mu_{E,V}  .  $$

\noindent 3) Let $A\in C(E,M^q)$; then, 
$$\frac{1}{\b^k_{E,V}}\L_{E,V}^k A\tends Q_{E,V}\inn{A}{\mu_{E,V}}   $$
uniformly as $k\tends+\infty$. If $A\in LogC^{a,\a}_+(E)$, then the convergence above is exponentially fast. A consequence of this is the last assertion of point 1), i.e. that the eigenvalue $\b_{E,V}$ of point 1) is simple. 

\noindent 4) $\mu_{E,V}(E)$ is a positive-definite matrix. 

\noindent 5) The Gibbs property holds; in other words, there is $D_1>0$ such that 
$$\frac{e^{V^l(x)}}{D_1\b^l_{E,V}}\tr\Psi_{x_0\dots x_{l-1}}\cdot\mu_{E,V}(E)\le
\mu_{E,V}([x_0\dots x_{l-1}]_{\a})\le
D_1\cdot\frac{e^{V^l(x)}}{\b^l_{E,V}}\tr\Psi_{x_0\dots x_{l-1}}\cdot\mu_{E,V}(E)  .  $$
In the formula above, $[x_0\dots x_{l-1}]_\a=[x_0\dots x_{l-1}]$ if we are on $\Sigma$ and 
$[x_0\dots x_{l-1}]_\a=[x_0\dots x_{l-1}]_G$ if we are on $G$. We have also set  
$$V^l(x)=V(x)+V\circ\s(x)+\dots+V\circ\s^{l-1}(x) $$
if we are on $\Sigma$ and 
$$V^l(x)=\tilde V(x)+\tilde V\circ F(x)+\dots+\tilde V\circ F^{l-1}(x) $$
if we are on $G$. Lastly, $\tr\Psi_{x_0\dots x_{l-1}}$ is defined as in (10). 

\noindent 6) Let us define $\kappa_{E,V}=(Q_{E,V},\mu_{E,V})_{HS}$; then, 
$\kappa_{E,V}$ is an atomless probability measure ergodic for $S$. Actually, the mixing property holds, which we state in the following way: if $g\in C(E,\R)$ and $A\in C(E,M^q)$, then  
$$\int_E(g\circ S^n\cdot A,\dr\mu_{E,V})_{HS}\tends
\int_E(gQ_{E,V},\dr\mu_{E,V})_{HS}\cdot
\int_E(A,\dr\mu_{E,V})_{HS}   .   $$
The convergence is exponentially fast if $A\in LogC^{a,\a}_+(E)$ and $g\in C^{0,\a}(E,\R)$. 

\rm

\vskip 1pc

\noindent{\bf Definition.} We shall say that $\mu_{E,V}$ is the Gibbs measure on $E$ and that $\kappa_{E,V}=(Q_{E,V},\mu_{E,V})_{HS}$ is Kusuoka's measure on $E$. Since 
$\mu_{E,V}\in\pc_{Q_{{E,V}}}$, $\kappa_{E,V}$ is a probability measure. 

\vskip 1pc

We recall Lemma 4.6 of [3], which shows that there is a natural relationship between Gibbs measures on $\Sigma$ and $G$; the notation is that of Proposition 1.2. 

\lem{1.3} We have that $\b_{G,V}=\b_{\Sigma,V}$; we shall call $\b_V$ their common value. Up to multiplying one of them by a positive constant, we have that 
$Q_{\Sigma,V}=Q_{G,V}\circ\Phi$; the conjugation $\fun{\Phi}{\Sigma}{G}$ is that of (1.11).  
For this choice of $Q_{\Sigma,V}$ and $Q_{G,V}$, let $\mu_{\Sigma,V}$ and $\mu_{G,V}$ be as in point 2) of Proposition 1.2; then $\mu_{G,V}=\Phi_\sharp\mu_{\Sigma,V}$ and 
$\kappa_{G,V}=\Phi_\sharp\kappa_{\Sigma,V}$, where $\Phi_\sharp$ denotes the push-forward by $\Phi$. 

\rm

\vskip 1pc

\noindent{\bf Definition.} Let $\b_V$ be the eigenvalue of the lemma above; we shall say that $P(V)\colon=\log\b_V$ is the pressure of the potential $V$.

\vskip 1pc

Just a word on the relationship with [5] and [13]. It is standard ([2]) that, if $V\equiv 0$, then the eigenfunction $Q$ is constant and $D_1=1$ in point 5) of Proposition 1.2. On constant matrices, our operator $\L_{\Sigma,V}$ and the operator $L_A$ of Proposition 15 of [13] coincide; in particular, they have the same spectral radius and the same maximal eigenvector, $Q$. Looking at the eigenvalues, this means that our $\b_{\Sigma,V}$ coincides with their $e^{P(A,2)}$. Moreover, when $V\equiv 0$, 
$\kappa_{\Sigma,V}$ coincides with the measure $m$ on the shift defined in [5] and [13]; we briefly prove why this is true.

A theorem of [5], recalled as Theorem 2 in [13], says that the probability measure $m$ of [13] satisfies the following: there are $C>0$ and $P(A,2)\in\R$ such that, for all 
$x_1,\dots,x_n\in\{ 1,\dots,t \}$, we have 
$$\frac{1}{C}m([x_1\dots x_n])\le
\frac{||A_{x_n}\cdot\dots\cdot A_{x_1}||_{HS}^2}{e^{nP(A,2)}}\le
Cm([x_1\dots x_n])  .  $$
We have already seen that $P(A,2)=\b_V$. Together with the Gibbs property of Proposition 1.2, the formula above implies that $m$ is absolutely continuous with respect to 
$\kappa_{\Sigma,V}$. Since $m$ is invariant and $\kappa_{G,V}$is ergodic, they must coincide. 

Before stating the last lemma of this section we need to define the Lyapounov exponents of the fractal. We consider the map 
$$\fun{H}{G}{M^d}$$
$$\fun{H}{x}{(D\psi_i)_\ast} \txt{if}x\in\oc_i    $$
which is defined $\kappa_{G,V}$-a.e. on $G$. For the expanding map $F$ defined in (F4) we set 
$$H^l(x)=H(x)\cdot H(F(x))\cdot\dots\cdot H(F^{l-1}(x))  .  $$
Since Kusuoka's measure $\kappa_{G,V}$ is ergodic for the map $F$, we can apply the results of [21] (see [8] for a more elementary presentation) and get that for 
$\kappa_{G,V}$-a.e. $x\in G$ there is $\Lambda_x\in M^d$ such that 
$$\lim_{l\tends+\infty}
[H^l(x)\cdot\mu_{G,V}(G)\cdot\tr H^l(x)]^\frac{1}{2l}=\Lambda_x  .  \eqno (1.25)$$
The fact that $\kappa_{G,V}$ is ergodic for $F$ implies by [21] that neither the eigenvalues of $\Lambda_x$ nor their multiplicities depend on $x$. 

Though we won't need it in the following, we recall Lemma 6.2 of [3]. 

\lem{1.4} Let the maps $\psi_i$, the Gibbs measure $\mu_{G,V}$ and Kusuoka's measure 
$\kappa_{G,V}$ be as in Theorem 1.1. Let us suppose that the maximal eigenvalue of 
$\Lambda_x$ is simple and let us call $\hat E_x$ its eigenspace; let $P_{\hat E_x}$ be the orthogonal projection on $\hat E_x$. Then, we have that 
$$\mu_{G,V}=P_{\hat E_x}||\mu_{G,V}||  .  $$
Moreover, we have that, for $\kappa_{G,V}$-a.e. $x\in G$, 
$$P_{\hat E_x}=\lim_{l\tends+\infty}
\frac{H^l(x)\cdot\mu_{G,V}(G)\cdot\tr H^l(x)}{|| H^l(x)\cdot\mu_{G,V}(G)\cdot\tr H^l(x) ||_{HS}}  .  
\eqno (1.26)$$

\rm

\vskip 2pc

\centerline{\bf \S 2}
\centerline{\bf Pressure}

\vskip 1pc

We define the map $\fun{\g_i}{\Sigma}{\Sigma}$ by $\g_i(x)=ix$. As in the last section, 
$(E,S)$ is either one of $(\Sigma,\s)$ or $(G,F)$; in both cases we denote by $\a_i$ a branch of the inverse of $S$, i.e. $\a_i=\g_i$ if $E=\Sigma$ and $\a_i=\psi_i$ if $E=G$. If $m$ is a probability measure on $E$ invariant for the map $S$, we denote by $h(m)$ its entropy. 

Given a potential $V\in C^{0,\a}$(E,\R), let $Q_{E,V}$ be as in Proposition 1.2;  we define 
$\kc_V(E)$ as the set of the couples $\{ m,M \}$ where $m$ is a non-atomic, $S$-invariant probability measure on $E$ and $\fun{M}{E}{M^q}$ is a Borel function such that 
$$M(x)\txt{ is positive semi-definite and}
(Q_{E,V}(x),M(x))_{HS}=1\txt{for $m$-a.e. $x\in E.$}  $$
As a consequence of the equality above, defining $\pc_{Q_{E,V}}$ as in (1.23) we have 
$$m=(Q_{E,V},Mm)_{HS}\txt{and thus}Mm\in \pc_{Q_{E,V}}   .  \eqno (2.1)$$
Let now $\kappa_{E,V}$ and $\mu_{E,V}$ be as in Proposition 1.2; recalling that 
$||\mu_{E,V}||$ and $\kappa_{E,V}$ are mutually absolutely continuous by (1.22), we can write 
$$\mu_{E,V}=\tilde M_{E,V}\cdot||\mu_{E,V}||=
M_{E,V}\cdot\kappa_{E,V}  $$
with $M_{E,V}=\tilde M_{E,V}\cdot\frac{\dr||\mu_{E,V}||}{\dr\kappa_{E,V}}$, where the derivative is in the Radon-Nikodym sense. This yields the first equality below, while the second one comes from the definition of $\kappa_{E,V}$ in point 6) of Proposition 1.2. 
$$(Q_{E,V},M_{E,V})_{HS}\cdot\kappa_{E,V}=
(Q_{E,V},\tilde M_{E,V})_{HS}\cdot||\mu_{E,V}||=\kappa_{E,V}  .  $$
As a consequence, we get that 
$$(Q_{E,V}(x),M_{E,V}(x))_{HS}= 1\txt{for $\kappa_{E,V}$-a.e. $x\in G$,}    $$
i.e. that $(\kappa_{E,V},M_{E,V})\in\kc_V(E)$.

We denote by $\bc$ the Borel $\s$-algebra of $E$ and by $m[\a_i(E)|S^{-1}(x)]$ the conditional expectation of $1_{\a_i(E)}$ with respect to the sigma-algebra $S^{-1}(\bc)$; paradoxically, this notation means that $m[\a_i(E)|S^{-1}(x)]$ is a Borel function of $S(x)$. 

We denote by $h(m)$ the entropy of a $S$-invariant measure on $E$. For 
$\{ m,M \}\in\kc_V(E)$ we set 
$$I_E(m,M)\colon=
h(m)+\int_EV(x)\dr m(x)+$$
$$\sum_{i=1}^t
\int_E\log(Q_{E,V}(x),\tr\Psi_i\circ M\circ S(x))_{HS}\dr m(x)  . \eqno (2.2)$$
By the last two formulas, (7) becomes 
$$I_G(\kappa_{G,V},M_{G,V})=\sup_{\{ m,M \}\in\kc_V(G)}I_G(m,M)  .  $$
The definition of $I_E$ in (2.2) implies immediately that $I_E$ is conjugation-invariant; in particular, if $\fun{\Phi}{\Sigma}{G}$ is the coding of Section 1 and 
$(m,M)\in\kc_V(\Sigma)$, then
$$I_\Sigma(m,M)=I_G(\Phi_\sharp m,M\circ\Phi^{-1})  ,  $$
where $\Phi_\sharp m$ denotes the push-forward; note that $\Phi^{-1}$ is defined 
$m$-a.e.: indeed, $(m,M)\in\kc_V(G)$ implies that $m$ is non-atomic and we have shown in Section 1 that the set $\{ y\in G\st\sharp\Phi^{-1}(y)>1 \}$ is countable. 

By the two formulas above and Lemma 1.3, (7) is equivalent to 
$$I_\Sigma(\kappa_{\Sigma,V},M_{\Sigma,V})=
\sup_{\{ m,M \}\in\kc_V(\Sigma)}I_\Sigma(m,M)    \eqno (2.3)$$
and this is what we shall prove. 

\vskip 1pc

\noindent{\bf Notation.} Since from now on we shall work on $\Sigma$, we simplify our notation:  
$$\L_V\colon=\L_{\Sigma,V},\qquad
\kappa_V\colon=\kappa_{\Sigma,V},\qquad
\mu_V\colon=\mu_{\Sigma,V},\qquad
\b_V\colon=\b_{\Sigma,V},\qquad
P_V\colon=\log\b_{\Sigma,V}, \qquad
Q_V\colon=Q_{\Sigma, V}  .  $$

\vskip 1pc

We begin writing $I_\Sigma$ in a slightly different way. Since $m[\g_i(E)|\s^{-1}(x)]$ is the conditional expectation of $1_{\g_i(\Sigma)}$ with respect to $\s^{-1}(\bc)$, we get the first equality below; the second one follows since $\g_i\circ\s$ is the identity on $\g_i(\Sigma)$; the last equality follows since the sets $\g_i(\Sigma)$ are a partition of $\Sigma$. 
$$\sum_{i=1}^t\int_\Sigma m[\g_i(\Sigma)|\s^{-1}(x)]\cdot\log(Q_{V}\circ\g_i\circ\s(x),\tr\Psi_i\cdot M\circ\s(x))_{HS}
\dr m(x)=$$
$$\sum_{i=1}^t\int_{\g_i(\Sigma)}\log(Q_{V}\circ\g_i\circ\s(x),\tr\Psi_i\circ M\circ\s(x))_{HS}\dr m(x)  =  $$
$$\sum_{i=1}^t
\int_{\g_i(\Sigma)}\log(Q_{V}(x),\tr\Psi_i\circ M\circ\s(x))_{HS}\dr m(x)=$$
$$\int_\Sigma\log(Q_{V}(x),\tr\Psi_i\circ M\circ\s(x))_{HS}\dr m(x)  .  $$
We can thus re-write (2.2) as 
$$I_\Sigma(m,M)\colon=
h(m)+\int_\Sigma V(x)\dr m(x)+$$
$$\sum_{i=1}^t\int_\Sigma m[\g_i(\Sigma)|\s^{-1}(x)]\cdot\log(Q_{V}\circ\g_i\circ\s(x),\tr\Psi_i\cdot M\circ\s(x))_{HS}
\dr m(x)  .  \eqno (2.4)$$

Let $m$ be a $\s$-invariant probability measure on $\Sigma$; in Formula (2.5) below we express $m[i|\s^{-1}(x)]$ using disintegration, in (2.6) we express it as a limit. We recall that a statement of the disintegration theorem is in Theorem 5.3.1 of [1]; the full theory is in [19], while [17] and [18] follow a different approach. 

We disintegrate $m$ with respect to the map $\fun{\s}{\Sigma}{\Sigma}$, which is measurable from $\bc$ to $\bc$ and preserves $m$. The disintegration theorem implies the first equality below, the second one follows from the fact that $m$ is 
$\s$-invariant. 
$$m=\tilde\eta_{x}\otimes(\s_\sharp m)=\tilde\eta_{x}\otimes m  .  $$
This means the following: first, $\tilde\eta_{x}$ is a probability measure on $\Sigma$ which concentrates on $\s^{-1}(x)$. Second, if $B\in\bc$, then the map $\fun{}{x}{\tilde\eta_{x}(B)}$ is $\bc$-measurable. Lastly, if $\fun{f}{\Sigma}{\R}$ is a bounded $\bc$-measurable function, then
$$\int_{\Sigma}f(x)\dr m(x)=
\int_\Sigma\dr\s_\sharp( m)(x)\int_\Sigma f(z)\dr\tilde\eta_{x}(z)  .  $$
We set $\eta_x=\tilde\eta_{\s(x)}$; by the remarks above, $\eta_x$ concentrates on 
$\s^{-1}(\s(x))$, i.e. on the fiber containing $x$, and, if $B\in\bc$, the map 
$\fun{}{x}{\eta_x(B)}$ is $\s^{-1}(\bc)$-measurable. 

Let $\fun{g}{\Sigma}{\R}$ be a bounded Borel function; the first equality below comes from the disintegration above; the second one follows from the fact that $\s_\sharp m=m$ and the definition of $\eta_x$. The third equality follows from the fact that $g\circ\s$ is constantly equal to $g\circ\s(x)$ on the fiber $\s^{-1}(\s(x))$ on which 
$\eta_x=\tilde\eta_{\s(x)}$ concentrates. 
$$\int_\Sigma 1_{[i]}(x) g\circ\s(x)\dr m(x)=
\int_\Sigma\dr m(x)\int_\Sigma 1_{[i]}(z)g\circ\s(z)\dr\tilde\eta_x(z)=$$
$$\int_\Sigma\dr m(x)\int_\Sigma 1_{[i]}(z)g\circ\s(z)\dr\eta_x(z)=
\int_\Sigma g\circ\s(x)\dr m(x)\int_\Sigma 1_{[i]}(z)\dr\eta_x(z)  .  $$
Since $g\circ\s$ is an arbitrary, bounded $\s^{-1}(\bc)$-measurable function and we saw above that $\eta_x([i])$ is $\s^{-1}(\bc)$-measurable, the formula above is the definition of conditional expectation. In other words, we have the first equality below, while the second one is trivial. 
$$m[i|\s^{-1}(x)]=\eta_{x}([i])=\int_\Sigma 1_{[i]}(z)\dr\eta_{x}(z)  .  \eqno (2.5)$$

Rokhlin's theorem ([17], [18]) is an analogue of Lebesgue's differentiation theorem for sequences of finer and finer partitions; the ones we consider are 
$$\{
\s^{-1}[x_0\dots x_l]\st \{ x_i \}_{i\ge 0}\in\Sigma\qquad l\ge 0
\}  .     $$
These partitions generate $\s^{-1}(\bc)$; by Rokhlin's theorem, for $m$-a.e. $x\in\Sigma$, $\eta_x([i])$ is the limit of the quotient of the measure of $[i]$ intersected with the element of the partition containing $\s^{-1}(\s(x))$ (i.e. the fiber containing $x$), and the measure of this element. This yields the first equality below, while the second one follows by the definition of $\s$ and the third one from the fact that $m$ is $\s$-invariant. 
$$\eta_{x}([i])=
\lim_{l\tends+\infty}
\frac{
m([i]\cap\s^{-1}(\s[x_0\dots x_l]))
}{
m(\s^{-1}(\s[x_0\dots x_l]))
}=$$
$$\lim_{l\tends+\infty}\frac{
m([ix_1\dots x_l])
}{
m(\s^{-1}([x_1\dots x_l]))
}  =
\lim_{l\tends+\infty}\frac{
m([ix_1\dots x_l])
}{
m([x_1\dots x_l])
}  .  $$
We already know it, but from the formula above it follows again that $\eta_x([i])$ is 
$\s^{-1}(\bc)$-measurable. Expressing $\eta_x([i])$ as a limit as in the formula above and using (2.5) we get that, for $m$-a.e. 
$x\in\Sigma$, 
$$m[i|\s^{-1}(x)]=
\lim_{l\tends+\infty}\frac{
m([ix_1\dots x_l])
}{
m([x_1\dots x_l])
}  .  \eqno (2.6)$$

We recall a lemma of [14]. 

\lem{2.1} Let $m$ be an invariant probability measure on $\Sigma$; then, for $m$-a.e. 
$x\in\Sigma$ we have 
$$\sum_{i=1}^t m[i|\s^{-1}(x)]=1  .  \eqno (2.7)$$

\proof The first equality below is (2.6) and it holds for $m$-a.e. $x\in\Sigma$; the second one follows since the cylinders $[ix_1\dots x_l]$ are disjoint as $i$ varies in $(1,\dots,t)$; the third one follows by the definition of the map $\s$ and the last one from the fact that 
$m$ is $\s$-invariant. 
$$\sum_{i=1}^t m[i|\s^{-1}(x)]=
\lim_{l\tends+\infty}\sum_{i=1}^t\frac{
m([ix_1\dots x_l])
}{
m([x_1\dots x_l])
}   =  $$
$$\lim_{l\tends+\infty}
\frac{
m\left(
\bigcup_{i=1}^t[ix_1\dots x_l]
\right)  
}{
m([x_1\dots x_l])
}
=
\lim_{l\tends+\infty}
\frac{
m(
\s^{-1}([x_1\dots x_l])}{
m([x_1\dots x_l])
} 
  =1  .  $$

\fin

We also recall from [14] a consequence of the Kolmogorov-Sinai formula ([24]): if $m$ is an invariant probability measure on $\Sigma$, then 
$$h(m)=
-\sum_{i=1}^t\int_\Sigma
m[i|\s^{-1}(x)]\cdot\log m[i|\s^{-1}(x)]\dr m(x)  .  \eqno (2.8)$$

The next lemma is another well-known fact. 

\lem{2.2} Let $m$ be an invariant probability measure on $\Sigma$; let 
$\fun{\g_i}{(x_0x_1\dots)}{(ix_0x_1\dots)}$ and let us suppose that 
$$(\g_i)_\sharp m=am   \eqno (2.9)$$
for a probability density $a$ which is positive on $[i]$. Then,  
$$m[i|\s^{-1}(x)]=\frac{1}{a\circ\g_i\circ\s(x)}  .  \eqno (2.10)$$

\proof The first equality below is (2.6) and it holds for $m$-a.e. 
$x=(x_0x_1\dots)\in\Sigma$; the second one follows from the definition of $\g_i$ and of the push-forward; the third one follows from (2.9) and the fourth one from Rokhlin's theorem applied to the partition of $[i]$ given by the sets $[ix_1\dots x_l]$ as $l\ge 1$ and 
$(ix_0x_1\dots)\in\Sigma$; it holds for $m$-a.e. $x\in \Sigma$. The last equality follows from the definition of $\g_i$ and $\s$. 
$$m[i|\s^{-1}(x)]=
\lim_{l\tends+\infty}\frac{
m([ix_1\dots x_l])
}{
m([x_1\dots x_l])
}  =
\lim_{l\tends+\infty}\frac{
m([ix_1\dots x_l])
}{
[(\g_i)_\sharp m]([ix_1\dots x_l])
}  =$$
$$\lim_{l\tends+\infty}
\frac{m([ix_1\dots x_l])}{\int_{[ix_1\dots x_l]}a(z)\dr m(z)}
=
\frac{1}{a(ix_1x_2\dots)}=
\frac{1}{a\circ\g_i\circ\s(x)}   .  $$

\fin

\lem{2.3} Let $V\in C^{0,\a}(\Sigma,\R)$, let $\mu_V$ and $Q_V$ be as in the notation after Formula (2.3); let $\fun{\g_i}{\Sigma}{\Sigma}$ be as in Lemma 2.2. Then, the following formulas hold. 
$$(\g_i)_\sharp\mu_V=\b_V e^{-V(x)}\cdot(\tr\Psi_i)^{-1}\cdot\mu_V|_{[i]}  .  
\eqno (2.11)$$
$$(\g_i)_\sharp\kappa_V=\b_V e^{-V(x)}(Q_V\circ\s(x),(\tr\Psi_i)^{-1}\cdot\mu_V|_{[i]})_{HS}  
\eqno (2.12)$$
or equivalently, setting $\mu_V=M_V\cdot\kappa_V$, 
$$(\g_i)_\sharp\kappa_V=1_{[i]}(x)\cdot \b_V e^{-V(x)}(Q_V\circ\s(x),(\tr\Psi_i)^{-1}\cdot M_V(x))_{HS}
\cdot\kappa_V  .   \eqno (2.13)$$
Lastly, 
$$\kappa_V[i|\s^{-1}(x)]=
\frac{
e^{V(i\s(x))}
}{
\b_V(Q_V\circ\s(x),(\tr\Psi_i)^{-1}\cdot M_V(i\s(x)))_{HS}
}  .  \eqno (2.14)$$

\proof We begin with (2.11). We recall that the Gibbs property follows iterating Lemma 4.1 of [3], which says the following: let $\eta\in(0,1)$ be as in (1.2); then, there is $D_1>0$ such that, for all $x=(x_0x_1\dots)\in\Sigma$ and $l\in\N$, 
$$\frac{1}{\b_V}e^{-D_1\eta^{\a l}}\cdot e^{V(x)}\cdot\tr\Psi_i\cdot\mu_V([x_1\dots x_{l-1}])\le
\mu_V([x_0x_1\dots x_{l-1}])\le$$
$$\frac{1}{\b_V}e^{D_1\eta^{\a l}}\cdot e^{V(x)}\cdot\tr\Psi_i\cdot\mu_V([x_1\dots x_{l-1}])  ,
\eqno (2.15)$$
where the inequalities are the standard ones between symmetric matrices. 

The first equality below is the definition of the push-forward, the second one follows by the definition of $\g_i$.
$$[(\g_i)_\sharp\mu_V]([x_0x_1\dots x_{l-1}])=\mu_V(\g_i^{-1}([x_0x_1\dots x_{l-1}]))
=\left\{
\eqalign{
\mu_V([x_1\dots x_{l-1}])  &\txt{if}x_0=i\cr
0 &\txt{otherwise.}  
}
\right.$$
Since both $(\g_i)_\sharp\mu_V$ and $\mu_V|_{[i]}$ concentrate on $[i]$, it suffices to show that the two measures of (2.11) coincide on the Borel sets of $[i]$. Thus, in the formula above we suppose that $x_0=i$ and we re-write (2.15) in this case.  
$$\frac{1}{\b_V}e^{-D_1\eta^{\a l}}\cdot e^{V(x)}\cdot\tr\Psi_i\cdot
[(\g_i)_\sharp\mu_V]([x_0x_1\dots x_{l-1}])\le
\mu_V([x_0x_1\dots x_{l-1}])\le$$
$$\frac{1}{\b_V}e^{D_1\eta^{\a l}}\cdot e^{V(x)}\cdot\tr\Psi_i\cdot
[(\g_i)_\sharp\mu_V]([x_0x_1\dots x_{l-1}]) . $$
This is true for all $x\in[x_0x_1\dots x_{l-1}]$; by the last formula and monotonicity of the integral we get that 
$$\frac{1}{\b_V}e^{-D_1\eta^{\a l}}
\int_{[x_0x_1\dots x_{l-1}]}e^{V(x)}\cdot\tr\Psi_i\cdot\dr[(\g_i)_\sharp\mu_V](x)\le
\mu_V([x_0\dots x_{l-1}])\le$$
$$\frac{1}{\b_V}e^{D_1\eta^{\a l}}
\int_{[x_0x_1\dots x_{l-1}]}e^{V(x)}\cdot\tr\Psi_i\cdot\dr[(\g_i)_\sharp\mu_V](x)  .   $$
Fix $k\ge 1$; since $[x_0\dots x_k]$ is the disjoint union of the cylinders 
$[x_0\dots x_k x_{k+1}\dots x_l]$, the last formula and the fact that 
$e^{D_1\eta^{\a l}}\tends 1$ as $l\tends+\infty$ imply easily that, if $x_0=i$ and $k\ge 1$, then 
$$\mu_V([x_0\dots x_{k-1}])=
\frac{1}{\b_V}\int_{[x_0x_1\dots x_{k-1}]}
e^{V(x)}\cdot\tr\Psi_i\cdot\dr[(\g_i)_\sharp\mu_V](x)  .  $$
Since the cylinders $[x_0\dots x_{k-1}]$ generate the Borel sets of $[x_0]=[i]$, the last formula implies that 
$$\mu_V|_{[i]}=\frac{1}{\b_V}e^V\cdot\tr\Psi_i\cdot[(\g_i)_\sharp\mu_V]$$
from which (2.11) follows. 

Now to (2.12). Let $\fun{f}{\Sigma}{\R}$ be a bounded Borel function. The first equality below is the definition of the push-forward, the second one of $\kappa_V$; the third one comes because, by definition, $\s\circ\g_i$ is the identity on $\Sigma$. The fourth one follows from the definition of the push-forward and the last one from (2.11). 
$$\int_{[i]} f(x)\dr[(\g_i)_\sharp\kappa_V](x)=
\int_{\Sigma} f\circ\g_i(z)\dr\kappa_V(z)=$$
$$\int_\Sigma f\circ\g_i(z)(Q_V(z),\dr\mu_V(z))_{HS}=
\int_\Sigma f\circ\g_i(z)(Q_V\circ\s\circ\g_i(z),\dr\mu_V(z))_{HS}=$$
$$\int_{[i]}f(x)\cdot
(Q_V\circ\s(x),\dr[(\g_i)_\sharp\mu_V](x))_{HS}=
\int_{[i]}f(x)\cdot\b_V e^{-V(x)}\cdot
(Q_V\circ\s(x),(\tr\Psi_i)^{-1}\cdot\dr\mu_V(x))_{HS}  .  $$

This shows (2.12); as for Formula (2.13), it follows from (2.12) and the fact that 
$\mu_V=M_V\cdot\kappa_V$. 

Formula (2.14) now follows from (2.13) and (2.10). 

\fin

We recall a well-known fact: if 
$$p_i,q_i\ge 0,\quad
\sum_{i=1}^t p_i=\l>0,\quad
\sum_{i=1}^t q_i=1  ,  \eqno (2.16)$$
then
$$-\sum_{i=1}^t q_i\log q_i+\sum_{i=1}^t q_i\log p_i\le\log\l  \eqno (2.17)$$
with equality if and only if $q_i=\frac{p_i}{\l}$ for all $i\in(1,\dots,t)$. 

Usually (see for instance [14] or [24]) this formula is stated for $\l=1$; applying  this particular case to $\left\{ \frac{p_i}{\l} \right\}_i$ we get the general one. 

For completeness' sake, we prove below another well-known formula from page 36 of [14]: if $\fun{g}{\Sigma}{\R}$ is a continuous function and $m$ is an invariant probability measure, then the second equality below comes from the definition of conditional expectation; the last one comes from the fact that the sets $[i]$ form a partition of 
$\Sigma$.
$$\sum_{i=1}^t\int_\Sigma g(ix_1\dots)m[i|\s^{-1}(x)]\dr m(x)=
\sum_{i=1}^t\int_\Sigma g\circ\g_i\circ\s(x) m[i|\s^{-1}(x)]\dr m(x)=$$
$$\sum_{i=1}^t\int_\Sigma g\circ\g_i\circ\s(x) 1_{[i]}(x)\dr m(x)=
\sum_{i=1}^t\int_{[i]}g(ix_1\dots)\dr m(x)=
\int_\Sigma g(x)\dr m(x)  .  \eqno (2.18)$$

\vskip 1pc

\noindent{\bf Proof of Theorem 1.} As we saw at the beginning of this section, it suffices to show Formula (2.3), and this is what we shall do.

For $\{ m,M \}\in\kc_V(\Sigma)$  and $x\in\Sigma$ we define the following two functions of $x$. 
$$q_i(x)=m[i|\s^{-1}(x)],\qquad
p_i(x)=e^{V(i\s(x))}(Q_V(i\s(x)),\tr\Psi_i\cdot M\circ\s(x))_{HS}  .  \eqno (2.19)$$
By (2.7), $\{ q_i \}_i$ satisfies (2.16) for $m$-a.e. $x\in\Sigma$; we show that, setting 
$\l=\b_V$, $\{ p_i \}_i$ too satisfies (2.16) $m$-a.e.. The first equality  below comes from the definition of $p_i$ in the formula above, for the second one we transpose; the third one follows by the definition of $\L_V=\L_{\Sigma,V}$ in (1.24); the fourth one comes from the fact that $\L_{V}Q_V=\b_V Q_V$ and the last equality follows from the fact that 
$(m,M)\in\kc_V(\Sigma)$. 
$$\sum_{i=1}^t p_i=\sum_{i=1}^te^{V(i\s(x))}\cdot
(Q_V(i\s(x)),\tr\Psi_i\cdot M\circ\s(x))_{HS}= $$
$$\sum_{i=1}^t e^{V(i\s(x))}(\Psi_i\cdot Q_V(i\s(x)),M\circ\s(x))_{HS}  =$$
$$((\L_{V}Q_V)(\s(x)),M\circ\s(x))_{HS}=
\b_V\cdot (Q_V\circ\s(x),M\circ\s(x))_{HS}=
\b_V    .   $$
Since (2.16) is satisfied, Formula (2.17) holds, implying that for $m$-a.e. $x\in G$ we have  
$$-\sum_{i=1}^tq_i(x)\log q_i(x)+
\sum_{i=1}^t q_i(x)[V(i\s(x))+\log(Q_V(i\s(x)),\tr\Psi_i\cdot M\circ\s(x))_{HS}]\le
\log\b_V  .   $$
Integrating against $m$ we get that 
$$-\sum_{i=1}^t\int_\Sigma q_i(x)\log q_i(x)\dr m(x)+$$
$$\sum_{i=1}^t\int_\Sigma
q_i(x)[V(i\s(x))+\log( Q_V(i\s(x)),\tr\Psi_i\cdot M\circ\s(x))_{HS}]\dr m(x)\le
\log\b_V  .  $$
By the first one of (2.19) and (2.8) we get that 
$$-\sum_{i=1}^t\int_\Sigma q_i(x)\log q_i(x)\dr m(x)=h(m)  .   $$
The first one of (2.19) and (2.18) imply that 
$$\sum_{i=1}^t\int_\Sigma 
q_i(x) V(i\s(x))\dr m(x)=
\int_\Sigma V(x)\dr m(x)  .  $$
The last three formulas imply that, if $(m,M)\in\kc_V(\Sigma)$ and $I_\Sigma$ is defined as in (2.4), then
$$I_\Sigma(m,M)\le\log\b_V  .  $$

To end the proof of (2.3)  it suffices to show that, in the formula above, equality is attained when $\{ m,M\}=\{ \kappa_V,M_V \}$. The only inequality in the derivation of the formula above is that of (2.17), which becomes an equality if  
$$q_i(x)=\frac{p_i(x)}{\b_V}\txt{for}i\in(1,\dots,t)\txt{and $\kappa_V$-a.e. }x\in\Sigma  .  
\eqno (2.20)$$
This is what we are going to prove. The first equality below is the definition of $q_i$ in (2.19) for $m=\kappa_V$, the second one is (2.14). 
$$q_i(x)=\kappa_V[i|\s^{-1}(x)]=
\frac{
e^{V(i\s(x))}
}{
\b_V(Q_V\circ\s(x),(\tr\Psi_i)^{-1}\cdot M_V(i\s(x)))_{HS}
}   .   \eqno (2.21)$$
We assert that (2.20) would follow if there were some positive Borel function 
$\fun{\l}{\Sigma}{\R}$ such that 
$$\tr\Psi_i\cdot M_V\circ\s(x)=\l(x) M_V(i\s(x))     \eqno (2.22)$$
for $m$-a.e. $x\in\Sigma$. 

Indeed, this formula implies the two outermost equalities below; the two innermost ones come from the fact that $\{ \kappa_V,M_V \} \in\kc_V(\Sigma)$. 
$$\frac{
1
}{
(Q_V\circ\s(x),\tr\Psi_i^{-1}\cdot M_V(i\s(x)))_{HS}
}  =
\frac{
\l(x)
}{
(Q_V\circ\s(x),M_V\circ\s(x))_{HS}
}  =\l(x)=$$
$$\l(x)(Q_V(i\s(x)),M_V(i\s(x)))_{HS}=
(Q_V(i\s(x)),\tr\Psi_i\cdot M_V\circ\s(x))_{HS}   .   $$
The first equality below is (2.21); the second one is the formula above; the last one is (2.19). 
$$q_i(x)=
\frac{
e^{V(i\s(x))}
}{
\b_V(Q_V\circ\s(x),(\tr\Psi_i)^{-1}\cdot M_V(i\s(x)))_{HS}
}=$$
$$\frac{1}{\b_V}e^{V(i\s(x))}\cdot 
(Q_V(i\s(x)),\tr\Psi_i\cdot M_V\circ\s(x))_{HS}=
\frac{p_i(x)}{\b_V}    $$
proving (2.20). 

In order to show (2.22), it suffices to recall that, for $\kappa_V$-a.e. $x\in\Sigma$, 
$$M_V(x)\cdot\frac{\dr\kappa_V}{\dr||\mu||}(x)
=\lim_{l\tends+\infty}
\frac{
\tr\Psi_{x_0\dots x_{l-1}}\cdot\mu_V(E)
}{
||\tr\Psi_{x_0\dots x_{l-1}}\cdot\mu_V(E)||_{HS}
} . $$
The proof of this equality follows from point 5) of Proposition 1.2 and Rokhlin's theorem; the details are in [3]. By (10), the last formula easily implies (2.22).

\fin

\vskip 2pc

\centerline{\bf \S 3}
\centerline{\bf The non weakly mixing case}

\vskip 1pc

In this section we give an estimate from above on the asymptotics of the periodic orbits of the shift when the potential $\hat V$ is constant. This is the so-called non-weakly mixing case of [14] and [15], in the sense that the shift is weakly mixing, but its suspension by 
$\hat V\equiv 1$ is not; we refer the reader to [14] for a precise definition of these objects. 

The first order of business is to show that counting periodic orbits on $\Sigma$ is the same as counting periodic orbits on $G$. We do this in Lemma 3.1 below, but we need a hypothesis stronger that (F3).

\vskip 1pc  

\noindent {\bf (F3+)}  We ask that, if  $i_0,\dots, i_{l-1},j_0,\dots,j_{k-1}\in(1,\dots,p)$, then 
$$\psi^{-1}_{i_0\dots i_{l-1}}\circ\psi_{j_0\dots j_{k-1}}\circ\psi_{i_0\dots i_{l-1}}(\fc)\cap\fc=\emptyset$$
unless $k=0$, when $\psi_{j_0\dots j_{k-1}}=id$ by definition and thus the composition of the three maps on the left is the identity. Note that, since $\psi_{j_0\dots j_{k-1}}$ is a diffeomorphism, (F3+) implies (F3) taking $l=0$. 

We briefly show that (F3+) holds for the harmonic Sierpinski gasket. Since 
$\psi_{i_0\dots i_{l-1}}$ is a diffeomorphism, we can re-write (F3+) as 
$$\psi_{j_0\dots j_{k-1}}\circ\psi_{i_0\dots i_{l-1}}(\fc)\cap
\psi_{i_0\dots i_{l-1}}(\fc)=\emptyset  .  \eqno (3.1)$$
We use the notation of figure 1 and distinguish two cases: in the first one, the cell 
$\psi_{j_0\dots j_{k-1}}\circ\psi_{i_0\dots i_{l-1}}(G)$ is not contained in the cell 
$\psi_{i_0\dots i_{l-1}}(G)$; then the two cells can intersect only at points of the type 
$\psi_{i_0\dots i_{l-1}}(R)$ with $R\in\{ A,B,C \}$; these points are not in 
$\psi_{i_0\dots i_{l-1}}(\fc)$ since $\psi_{i_0\dots i_{l-1}}$ is a diffeomorphism and 
$\{ A,B,C \}$ and $\fc$ are disjoint. 

The second case is when 
$\psi_{j_0\dots j_{k-1}}\circ\psi_{i_0\dots i_{l-1}}(G)\subset\psi_{i_0\dots i_{l-1}}(G)$; if 
$k\ge l$ this implies that $j_0=i_0,\dots,j_{l-1}=i_{l-1}$. Applying 
$\psi_{i_0\dots i_{l-1}}^{-1}$, (3.1) becomes 
$$\psi_{j_l\dots j_{k-1}}\circ\psi_{i_0\dots i_{l-1}}(\fc)\cap\fc=\emptyset  $$
which is true since (F3) holds on the harmonic gasket. If $k<l$, we note the following: if 
$\psi_{j_0\dots j_{k-1}}\circ\psi_{i_0\dots i_{l-1}}(G)\subset\psi_{i_0\dots i_{l-1}}(G)$, then the first $l$ indices of $j_0\dots j_{k-1}i_0\dots i_{l-1}$ must be $i_0\dots i_{l-1}$. In other words, we have that $j_0=i_0,\dots, j_{k-1}=i_{k-1}$ and also $i_0=i_k,\dots, i_{l-k-1}=i_{l-1}$. Now we apply $\psi_{i_0\dots i_{l-1}}^{-1}$ and (3.1) becomes 
$$\psi_{i_{l-k}\dots i_{l-1}}(\fc)\cap\fc=\emptyset$$
which is implied by (F3).

\lem{3.1} Let (F1)-(F2)-(F3+)-(F4) hold. Then, the map $\fun{\Phi}{\Sigma}{G}$ of Section 1 is a bijection between the periodic points of $(\Sigma,\s)$ and those of $(G,F)$. 

\proof We consider the set $H$ of (1.12); since after this formula we saw that the conjugation is injective on $\Phi^{-1}(H^c)$, it suffices to show the following:

\noindent 1) every periodic orbit of $F$ has a periodic coding, which is trivial, and

\noindent 2) the points of $H$ do not have a periodic coding (and thus, in particular, are not periodic.) 

For point 2), we begin to show that the points with a periodic coding are not in $\fc$. Indeed, if $x=(x_0x_1\dots)$ satisfies $x_i=x_{i+n}$ for all $i$, then by the second formula of (1.14) we get that $\psi_{x_0\dots x_{n-1}}\circ\Phi(x)=\Phi(x)$; now (F3) implies that 
$\Phi(x)\not\in \fc$.

But points with a periodic coding are not even in $\psi_{w_0\dots w_l}(\fc)$. Indeed, let $x=(x_0x_1\dots)$ be $n$-periodic and let us suppose that 
$\Phi(x)=\psi_{w_0\dots w_l}(z)$ for some $z\in\fc$; then, again by (1.14), 
$$\psi_{x_0\dots x_{n-1}}\circ\psi_{w_0\dots w_l}(z)=\psi_{w_0\dots w_l}(z)  .  $$
Applying to both sides the inverse of $\psi_{w_0\dots w_l}$ (which is a diffeomorphism) (F3+) shows that $z\not\in\fc$, as we wanted. 

\fin

Thus, from now on we shall count periodic orbits on $\Sigma$; since we work on $\Sigma$, we stick to the notation after (2.3), but often writing $\b$ instead of $\b_V$.

\vskip 1pc

\noindent{\bf Definition.} We say that two periodic orbits $x,x^\prime\in\Sigma$ are equivalent if there is $k\in\N$ such that $x_i^\prime=x_{i+k}$ for all $i\ge 0$. We call 
$\tau=\{ x \}$ the equivalence class of $x$; if $x$ has minimal period $n$, then $\tau$ has $n$ elements. Said differently, we are going to count the periodic orbits $\tau$, not the periodic points $( x_i )\in\Sigma$. 

\vskip 1pc

Let now $\fun{\Psi_i}{M^q}{M^q}$  be as in (2), let $x=(x_0x_1\dots)\in\Sigma$ and let 
$l\in\N$; we define $\tr\Psi_{x,l}$ as in (10). The definition of the Hilbert-Schmidt product and (2) easily imply that 
$$\tr\Psi_i(A)=(D\psi_i)_\ast\cdot A\cdot\tr(D\psi_i)_\ast  .  $$
Iterating this formula and recalling (10) we get that
$$\tr\Psi_{x,l}(A)=(D\psi_{x_0\dots x_{l-1}})_\ast\cdot A\cdot
\tr(D\psi_{x_0\dots x_{l-1}})_\ast   .  \eqno (3.2)$$

\lem{3.2} Let the operators $\fun{\tr\Psi_{x,l}}{M^q}{M^q}$ be defined as in (10) or, which is equivalent, (3.2); recall that $m={\rm dim}M^q$. Let $x\in\Sigma$ be $n$-periodic, with $n$ not necessarily minimal. Then, the linear operator $\fun{\tr\Psi_{x,l}}{M^q}{M^q}$ satisfies the five points below. 

\noindent 1) For all $l\ge 0$ we can order the eigenvalues $\a_{x,1},\dots,\a_{x,m}$ of 
$\tr\Psi_{x,l}$ in such a way that 
$$1>\a_{x,1}\ge |\a_{x,2}| \ge\dots\ge |\a_{x,m}| > 0  .  \eqno (3.3)$$ 
Note that the largest eigenvalue is real. 

\noindent 2)  In the same hypotheses, $\a_{x,i}=\a_{\s^j(x),i}$ for all $j\ge 0$ and all 
$i\in(1,\dots,m)$. Since $\a_{x,i}$ does not depend on the particular point of the orbit on which we calculate it, if $n$ is the minimal period and $\tau=\{ x \}$, then we can set 
$\a_{\tau,i}\colon=\a_{x,i}$. 

\noindent 3) Let $x\in\Sigma$ be $n$-periodic with $n$ minimal and let $k\in\N$. Then, the eigenvalues of $\tr\Psi_{x,kn}$ are 
$$\a_{x,1}^k,\a_{x,2}^k,\dots,\a_{x,m}^k   
\txt{with}|\a_{x,i}^k|\in(0,1)\quad\forall i.  $$
By point 2), they are all shift-invariant.  

\noindent 4) Let $x$ be a periodic orbit of period $n$ and let $k,j\in\N$; then, 
$${\rm tr}\tr\Psi_{x,kn}={\rm tr}\tr\Psi_{\s^j(x),kn}  .  $$
In other words, the trace of $\tr\Psi_{x,kn}$ is the same for all $x$ on the same orbit. 

\noindent 5) For a $n$-periodic $x\in\Sigma$ and $k\in\N$ we have that 
$$0\le{\rm tr}\tr\Psi_{x,kn}  .  $$

\proof We begin with Formula (3.3) of point 1): since by (3.2) $\tr\Psi_{x,l}$ is the composition of bijective operators, it is bijective too; thus, all its eigenvalues are non zero and we get the inequality on the right of (3.3). We shall show at the end that 
$\a_{x,1}$, the eigenvalue with largest modulus, is positive and smaller than 1. 

Point 2) follows because the eigenvalues of a product of invertible matrices do not change under cyclical permutations of the factors, which implies that the eigenvalues of 
$\tr\Psi_{x,n}=\tr\Psi_{x_0}\circ\dots\circ\tr\Psi_{x_{n-1}}$ coincide with those of 
$\tr\Psi_{\s(x),n}=\tr\Psi_{x_1}\circ\dots\circ\tr\Psi_{x_{n-1}}\circ\tr\Psi_{x_{0}}$.  

Point 3) follows from two facts: the first one is that, by (10) or (3.2), 
$\tr\Psi_{x,kn}=(\tr\Psi_{x,n})^k$; the second one is that the eigenvalues of $A^k$ are the $k$-th powers of the eigenvalues of $A$. Point 2) implies the last assertion about shift-invariance. 
  
Point 4) is an immediate consequence of point 3).

Next, to point 5).  If $A\in M^q$ is semi-positive definite, then $\tr\Psi_{x,l}(A)$ is semi-positive definite by (3.2). By (1.19) this implies that  
$$(\tr\Psi_{x,l}(A),A)_{HS}\ge 0  .  $$
If the matrix  $\tr\Psi_{x,l}$ on $M^q$ were self-adjoint, this would end the proof, but this is not the case. However, the last formula implies that ${\rm tr}\tr\Psi_{x,kn}\ge 0$ follows if we find an orthonormal base of $M^q$ all whose elements are semipositive-definite. A standard one is this: if $\{ e_i \}$ is an orthonormal base of $\Lambda^q(\R^d)$, just take 
$\2(e_i\otimes e_j+e_j\otimes e_i)$. 

We show that the maximal eigenvalue in (3.3) is a positive real. Let us consider the cone $C\subset M^q$ of positive-definite matrices; in the last paragraph we saw that 
$\tr\Psi_{x,l}(C)\subset C$. We would like to apply Theorem 1.1 to this situation, but we don't know whether the constant $D$ in its statement is finite. It is easy to see that, if $\e>0$, then $D<+\infty$ for the operator $\tr\Psi_{x,l}+\e Id$; if we apply Theorem 1.1 and let $\e\searrow 0$, we get that the eigenvalue $\a_{x,1}$ of maximal modulus of $\tr\Psi_{x,l}$ is a non-negative real, as we wanted. 

Lastly, we end the proof of (3.3) showing that all eigenvalues have modulus smaller than 1.  We shall use a version of the $\sup$ norm on 
$\Lambda^q(\R^d)$: if $\o\in\Lambda^q(\R^d)$, we set 
$$||\o||_{\sup}=
\sup\{
\o(v_1,\dots,v_q)\st ||v_1||\cdot\dots\cdot||v_q||\le 1
\}   .  $$
Since $D\psi_i$ is a contraction, (1) implies immediately that $(D\psi_i)_\ast$ contracts the $\sup$ norm: 
$$||(D\psi_i)_\ast\o||_{\sup}< ||\o||_{\sup}  .  $$
Since the spectrum of $\tr\Psi_i$ coincides with that of $\Psi_i$, it suffices to show that the modulus of the eigenvalues of $\Psi_i$ is smaller than 1, and this is what we shall do. Let 
$A\in M^q$ be such that $\Psi_iA=\l A$; $|\l|<1$ follows if we show that, for some norm 
$|| \cdot ||_{1}$ over $M^q$, we have 
$$||\Psi_iA||_{1}<||A||_1  .   $$
The function $||\cdot||_1$ defined below iss a norm on the space $M^q$ of symmetric matrices.   
$$||A||_1=\sup_{||\o||_{\sup}\le 1}|(A\o,\o)|  .  $$

The first equality below comes from the definition of the norm in the formula above and the definition of $\Psi_i$ in (2); the second one is the definition of the adjoint; the inequality comes from the fact that $(D\psi_i)_\ast$ contracts the norm: the set of the forms 
$(D\psi_i)_\ast \o$ with $||\o||_{\sup}\le 1$ is contained in a sphere of radius strictly smaller than 1. The last equality is again the definition of the norm.
$$||\Psi_i(A)||_1=
\sup_{||\o||_{\sup}\le 1}|(\tr(D\psi_i)_\ast A(D\psi_i)_\ast \o,\o)|=$$
$$\sup_{||\o||_{\sup}\le 1}|( A(D\psi_i)_\ast \o,(D\psi_i)_\ast \o)|<
\sup_{||\o||_{\sup}\le 1}|(A\o,\o)|=||A||_1  .  $$

\fin

From now on the argument is that of [14], with some minor notational changes.  

We begin to assume that the potential $V$ depends only on the first two coordinates, i.e. that 
$$V(x)=V(x_0,x_1)  .  \eqno (3.4)$$
We are going to see that in this case the maximal eigenvalue of the Ruelle operator coincides with that of a finite-dimensional matrix; another consequence will be that the zeta function of Formula (3.11) below is meromorphic.

We consider $(M^q)^t$, the space of column vectors
$$M=\left(
\matrix{
M_1\cr
\dots\cr
M_t
}
\right)   \eqno (3.5)$$
where the entries are matrices $M_i\in M^q$; we denote by $\hat d$ the dimension of  
$(M^q)^t$; clearly, $\hat d=m\cdot t$ for $m$ as in Lemma 3.2. For 
$j\in (1,\dots,t)$ we define 
$$\fun{L_V^j}{(M^q)^t}{M^q},\qquad
L^j_VM=\sum_{i=1}^t \Psi_i(M_i)e^{V(i,j)}  .  $$
With the notation of (3.5) and of the last formula we define an operator  
$$\fun{L_V}{(M^q)^t}{(M^q)^t}, \qquad
L_VM=\left(
\matrix{
L^1_VM\cr
\dots\cr
L^t_VM
}
\right)   .  $$
It is easy to decompose the matrix $L_V$ into $t\times t$ blocks from $M^q$ into itself; the $ji$-th block of $L_V$ is given by 
$$L_V^{ij}=e^{V(i,j)} \Psi_i  .  \eqno (3.6)$$
We define the positive cone 
$$(M^q)^t_+\subset (M^q)^t  $$
as the set of the arrays (3.5) such that $M_i$ is positive-definite for all 
$i\in(1,\dots,t)$. 

On the cone $(M^q)_+^t$ we define the hyperbolic distance $\theta$ as in Section 1. The nondegeracy hypothesis $(ND)_q$ of Section 1 implies easily that 
$${\rm diam}_\theta[L_V((M^q)^t_+)]<+\infty  .  $$
We skip the standard proof of the formula above and of the fact that $((M^q)^t_+,\theta)$ is complete; by (3.9) below, they are a consequence of the corresponding facts on  $LogC^{\a,a}_+(\Sigma)$, which are proven in [2] and [3]. 

By Theorem 1.1, these two facts imply that $L_V$ has a couple eigenvalue-eigenvector 
$$(\b^\prime,Q^\prime)\in(0,+\infty)\times(M^q)^t_+  .  \eqno (3.7)$$
An element $M\in (M^q)^t$ defines a function $M\in C(\Sigma,\R)$ which depends only on the first coordinate of $x=(x_0x_1\dots)$; we continue to denote it by $M$: 
$$M(x)\colon=M_{x_0}  .  \eqno (3.8)$$
Let $V$ be as in (3.4) and let $M$ depend only on the first coordinate as in (3.8); by (1.24) and (3.6) we easily see that, for all $x=(x_0x_1\dots)\in\Sigma$, 
$$(\L_VM)(x)=(L_VM)_{x_0}  .  \eqno (3.9)$$
As a consequence, $(\b^\prime,Q^\prime)$ of (3.7) is a couple eigenvalue-eigenvector for 
$\L_V$ too; by the uniqueness of Proposition 1.2 we get that, if we multiply $Q^\prime$ by a suitable positive constant, 
$$(\b_V,Q_V)=(\b^\prime,Q^\prime)  .  $$
In particular, $Q_V(x)=Q_{x_0}^\prime$, i.e. $Q_V$ depends only on the first coordinate. 

Recall that we have set $\hat d={\rm dim}(M^q)^t$; an immediate consequence of the Perron-Frobenius theorem is that all other eigenvalues of $L_V$ have modulus strictly smaller that $\b=e^{P(V)}$. Let's call them names:  
$$\l_{2,V},\l_{3,V},\dots,\l_{\hat d,V} \txt{with}
|\l_{j,V}|<\b=e^{P(V)} \txt{for}j\ge 2  .  \eqno (3.10)$$

\vskip 1pc

\noindent{\bf Definitions.} We define $Fix_n$ as the set of $(x_0x_1\dots)\in\Sigma$ such that $x_i=x_{i+n}$ for all $i\ge 0$; note that $n$ is any period, not necessarily the minimal one. 

We define the zeta function $\zeta(z,V)$ as 
$$\zeta(z,V)=
\exp\sum_{n\ge 1}\frac{z^n}{n}\sum_{x\in Fix_n}e^{V^n(x)}{\rm tr}\tr\Psi_{x,n}   
\eqno (3.11)$$
where $\tr\Psi_{x,n}$ is defined as in (3.2) and $V^n$ as in point 5 of Proposition 1.2. 

We recall Proposition 5.1 of [14]. 

\prop{3.3} Let $V\in C^{0,\a}(\Sigma,\R)$ and let $\b=e^{P(V)}$ be as in Proposition 1.2; then, the radius of convergence of the power series in (3.11) is 
$$R=e^{-P(V)}=\frac{1}{\b}  .  \eqno (3.12)$$

\proof We begin to note that the innermost sum in (3.11) is non-negative: indeed, all the summands are non-negative by point 5) of Lemma 3.2. Thus, we can take logarithms and (3.12) is implied by the first equality below, while the second one is the definition of $P(V)$. 
$$\lim_{n\tends+\infty}
\frac{1}{n}\log\sum_{x\in Fix_n}e^{V^n(x)}{\rm tr}\tr\Psi_{x,n}=\log\b=P(V)  .  \eqno (3.13)$$
We show this. Since $V\in C^{0,\a}(\Sigma,\R)$, for any given $\e>0$ we can find a function $\fun{W}{\Sigma}{\R}$ depending only on finitely many coordinates such that 
$$||W-V||_{\sup}<\e  .  \eqno (3.14)$$
Up to enlarging our symbol space, we can suppose that $W$ depends only on the first two coordinates, i.e. that it satisfies (3.4). Thus, we can define an operator $L_W$ as in (3.6); this formula implies easily that the $i$-th block on the diagonal of $L^n_W$ is given by 
$$\sum_{\{ x\in Fix_n,\quad x_0=i \}}e^{W^n(x)}\Psi_{x,n}    $$
where $x=(x_0x_1\dots)$, $\Psi_{x,n}=\Psi_{x_{n-1}}\circ\dots\circ\Psi_{x_0}$ and $W^n$ is defined as in Proposition 1.2. Since ${\rm tr}(\tr B)={\rm tr}(B)$, the first equality below follows; the second one comes from the formula above and the linearity of the trace; the last one follows from (3.10) for $L_W$. 
$$\sum_{x\in Fix_n}e^{W^n(x)}{\rm tr}\tr\Psi_{x,n}=
\sum_{x\in Fix_n}e^{W^n(x)}{\rm tr}\Psi_{x,n}=
{\rm tr}(L^n_W)=e^{nP(W)}+\l^n_{2,W}+\dots+\l^n_{\hat d,W}  .  \eqno (3.15)$$
Since we saw above that the left hand side is positive, we can take logarithms in the last formula; we get the equality below, while the limit follows from the inequality of (3.10) applied to $W$. 
$$\frac{1}{n}\log\sum_{x\in Fix_n}e^{W^n(x)}{\rm tr}\tr\Psi_{x,n}=
\frac{1}{n}\log[e^{nP(W)}+\l^n_{2,W}+\dots+\l^n_{\hat d,W}]\tends P(W)  .    $$
This is (3.13) for the potential $W$; it implies (3.13) for the potential $V$ if we show that, as $W\tends V$ in $C(\Sigma,\R)$, we have that 
$$P(W)\tends P(V)  \eqno (3.16)$$
and 
$$\sup_{n\ge 1}
\left\vert
\frac{1}{n}\log\sum_{x\in Fix_n}e^{W^n(x)}{\rm tr}\tr\Psi_{x,n}-
\frac{1}{n}\log\sum_{x\in Fix_n}e^{V^n(x)}{\rm tr}\tr\Psi_{x,n}
\right\vert\tends 0  .   \eqno (3.17)$$
Formula (3.16) follows immediately from the theorem of contractions depending on a parameter applied to the operator $\L_V$ on the cone $LogC^{\a,a}_+(\Sigma)$ with the hyperbolic metric. 

For (3.17) we begin to note that, by (3.14), 
$$e^{-n\e}\cdot e^{W^n(x)}\le
e^{V^n(x)}\le
e^{n\e}\cdot e^{W^n(x)}   .  $$
Since ${\rm tr}\tr\Psi_{x,n}\ge 0$ by Lemma 3.2, the formula above implies that 
$$e^{-n\e}\sum_{x\in Fix_n}e^{W^n(x)}{\rm tr}\tr\Psi_{x,n}\le
\sum_{x\in Fix_n}e^{V^n(x)}{\rm tr}\tr\Psi_{x,n}\le
e^{n\e}\sum_{x\in Fix_n}e^{W^n(x)}{\rm tr}\tr\Psi_{x,n}  .  $$
Taking logarithms, (3.17) follows and we are done. 

\fin

Now we suppose that (3.4) holds. The first equality below is the definition (3.11); the second and third ones come from formula (3.15) applied to the potential $V$. The next two equalities come from the properties of the logarithm, the last one comes from the definition of $\l_{i,V}$ in (3.10); the matrix $Id$ is the identity on $(M^q)^t$. 
$$\zeta(z,V)=
\exp\sum_{n\ge 1}\frac{z^n}{n}\sum_{x\in Fix_n}e^{V^n(x)}{\rm tr}\tr\Psi_{x,n}=$$
$$\exp\sum_{n\ge 1}\frac{z^n}{n}{\rm tr}(L^n_V)=
\exp\sum_{n\ge 1}\frac{z^n}{n}\left[
e^{nP(V)}+\l^n_{2,V}+\dots+\l^n_{\hat d,V}
\right]  =  $$
$$\exp[
-\log(1-ze^{P(V)})-\log(1-\l_{2,V}z)-\dots-\log(1-\l_{\hat d,V}z)
]=$$
$$\frac{
1
}{
(1-ze^{P(V)})(1-\l_{2,V}z)\cdot\dots\cdot(1-\l_{\hat d,V}z)
}   =  $$
$$\frac{1}{
{\rm det}(Id-zL_V)
}   .  \eqno (3.18)$$
Thus, if (3.4) holds (or if $V$ depends on finitely many coordinates, which we have seen to be equivalent), the function $\zeta(z,V)$ is meromorphic. 

\vskip 1pc

\noindent{\bf Definitions.} Let $x\in\Sigma$ be periodic of minimal period $n$ and let 
$\tau=\{ x \}$ be its equivalence class; with the notation of Lemma 3.2 we define
$$\l(\tau)\colon=\sharp\tau=n,\qquad
V(\tau)\colon=V^n(x),\qquad
\a_{\tau,i}\colon=\a_{x,i},\qquad
{\rm tr}\tr\Psi_\tau\colon={\rm tr}\tr\Psi_{x,n}=\a_{\tau,1}+\dots+\a_{\tau,m}  .  
\eqno (3.19)$$
Periodicity implies immediately that $V(\tau)$ does not depend on the representative 
$x\in\tau$ we choose.  As for $\a_{\tau,i}$ and ${\rm tr}\tr\Psi_\tau$, they too don't depend on the representative $x\in\tau$ by points 2) and 4) of Lemma 3.2.  

As a further bit of notation, we call $\tau^\prime$ a $k$-multiple of a prime orbit $\tau$. More precisely, $\tau^\prime=(\tau,k)$ where $\tau=\{ x \}$ is a periodic orbit of minimal period $n$ and $k\in\N$; we set 
$$\left\{
\eqalign{
\Lambda(\tau^\prime)&\colon=\l(\tau)=n\cr
\l(\tau^\prime)&\colon=kn\cr
{\rm tr}\tr\Psi_{\tau^\prime}&\colon={\rm tr}\tr\Psi_{x,kn}=\a^k_{\tau,1}+\dots+\a^k_{\tau,m}\cr
V(\tau^\prime)&\colon=V^{kn}(x)=kV(\tau)      .     
}
\right.   \eqno (3.20)$$
In the third formula above, the second equality is point 3) of Lemma 3.2; the notation for the eigenvalues is that of (3.19). In the last formula above, the second equality is an immediate consequence of periodicity. We stress again that none of the expressions above depends on the particular $x\in\tau$ we choose.  

Now we can give an alternate expression for the zeta function $\zeta(z,V)$; this time, we don't require (3.4). The first equality below is the definition of $\zeta(z,V)$ in (3.11); for the second one we group together in the innermost sum all the multiples of a periodic orbit. In the third one the outer sum is over the prime orbits $\tau$; we recall that $\tau$ has $n$ elements and use (3.20). The last two equalities follow from the properties of the logarithm. 
$$\zeta(z,V)=
\exp\sum_{n\ge 1}\frac{z^n}{n}\sum_{x\in Fix_n}e^{V^n(x)}{\rm tr} \tr\Psi_{x,n}=$$
$$\exp\sum_{n\ge 1}\frac{1}{n}\sum_{\{ \s^n(x)=x,\quad n\quad least \}}
\sum_{k\ge 1}e^{kV^n(x)}{\rm tr}(\tr\Psi_{x,kn})\cdot\frac{z^{nk}}{k}=$$
$$\exp\sum_\tau\sum_{k\ge 1}e^{kV(\tau)}[\a^k_{\tau,1}+\dots+\a^k_{\tau,m}]
\frac{z^{\l(\tau)k}}{k}=$$
$$\exp \left\{
-\sum_\tau\sum_{j=1}^{m}\log[1-\a_{\tau,j}e^{V(\tau)}\cdot z^{\l(\tau)}]
\right\}  =  $$
$$\prod_\tau\prod_{j=1}^{m}
\frac{
1
}{
1-\a_{\tau,j}e^{V(\tau)}\cdot z^{\l(\tau)}
}   .   $$
We single out two equalities from the formula above. 
$$\zeta(z,V)=\prod_\tau\prod_{j=1}^{m}
\frac{
1
}{
1-\a_{\tau,j}e^{V(\tau)}\cdot z^{\l(\tau)}
}  =  $$
$$\exp\sum_\tau\sum_{j=1}^{m}\sum_{k\ge 1}
\frac{\a^k_{\tau,j}\cdot e^{kV(\tau)}\cdot z^{k\l(\tau)}}{k}   .   \eqno (3.21)$$
By Proposition 3.3 we know that the series above converges if $|z|<\frac{1}{\b}$; by the way, this justifies the formal manipulations above. 

\noindent{\bf Proof of point 1) of theorem 2.} The first equality below comes from (3.18) and requires (3.4), for the second one we use (3.10). 
$$\zeta(z,V)=
\frac{
1
}{
{\rm det}(Id-zL_V)
}  =  $$
$$\frac{1}{1-\b z}\cdot
\prod_{j=2}^{\hat d}\frac{1}{1-\l_{j,V}z}  .  $$
In other words, 
$$\zeta(z,V)\cdot (1-\b z)=\tilde\a(z)      \eqno (3.22)$$
and, by the inequality of (3.10),  there is $\e>0$ such that $\tilde\a$ is non-zero and analytic in 
$$B_\e\colon=\left\{
z\in\C\st |z|<\frac{e^\e}{\b}
\right\}  .    \eqno (3.23)$$

The first equality below follows taking derivatives in the last expression of (3.21); recall that  $\l(\tau)$ is defined in (3.19). Differentiating the product in (3.22) and recalling that $\tilde\a\not=0$ in $B_\e$ we see that the second equality below holds for a function $\a$ analytic in the ball $B_\e$ of (3.23). 
$$\sum_\tau\sum_{j=1}^{m}\sum_{k\ge 1}\l(\tau)\cdot\a^k_{\tau,j}\cdot e^{kV(\tau)}\cdot
z^{k\cdot\l(\tau)-1}=\frac{
\zeta^\prime(z,V)
}{
\zeta(z,V)
}  =$$
$$\frac{\b}{1-\b z}+\a(z)  .  \eqno (3.24)$$

From now on we shall suppose that $V\equiv 0$. We re-write (3.24) under this assumption; for the term on the left we sum on the periodic orbits $\tau^\prime=(\tau,k)$ and use (3.20), for the term on the right we use the geometric series. 
$$\frac{1}{z}\sum_{\tau^\prime}\Lambda(\tau^\prime)
{\rm tr}\tr\Psi_{\tau^\prime}
\cdot z^{\l(\tau^\prime)}=
\frac{
\zeta^\prime(z,V)
}{
\zeta(z,V)
}  =$$
$$\frac{1}{z}\sum_{n\ge 1}(\b z)^n+\a(z)  .  $$
Subtracting the sum on the right, we get that  
$$\sum_{n\ge 1}z^{n-1}\left[
\sum_{\l(\tau^\prime)=n}\Lambda(\tau^\prime){\rm tr}\tr\Psi_{\tau^\prime}
-\b^n
\right]  =\a(z)  .  $$
Since $\a$ is analytic in the disc $B_\e$ of (3.23) we get that the radius of convergence of the series on the left is $\frac{e^\e}{\b}>\frac{1}{\b}$. Since the general term of a converging series tends to zero, we deduce that, for any smaller $\e>0$, 
$$\frac{e^{n\e}}{\b^n}\left[
\sum_{\l(\tau^\prime)=n}\Lambda(\tau^\prime){\rm tr}
\tr\Psi_{\tau^\prime}-\b^n
\right]  \tends 0\txt{as}n\tends+\infty  .  \eqno (3.25)$$
In particular, the sequence above is bounded by some $D_1>0$. 

We consider the case $\b>1$; we define $\eta(r)$ in the first equality below, where the sum is over the periodic orbits $\tau^\prime=(\tau,k)$, not just the prime ones; $\l$ and 
$\Lambda$ are defined in (3.20). The second equality comes (for $r$ integer) by adding and subtracting.
$$\eta(r)\colon=\sum_{\lambda(\tau^\prime)\le r}
\Lambda(\tau^\prime){\rm tr}\tr\Psi_{\tau^\prime}=$$
$$\sum_{n\in[1,r]}\left[
\sum_{\l(\tau^\prime)=n}\left(
\Lambda(\tau^\prime){\rm tr}\tr\Psi_{\tau^\prime}
\right)
-\b^n
\right]  +
\b\frac{\b^r-1}{\b-1}  .  \eqno (3.26)$$
By (3.26) and the remark after (3.25) there is $D_1>0$, independent of $r$, such that the first inequality below holds; the second one holds for some $D_2=D_2(\b)>0$ independent of $r$ by the properties of the geometric series; note that $\frac{\b}{e^\e}>1$ if $\e$ is small enough, since we are supposing that $\b>1$. 
$$\left\vert
\eta(r)-\frac{\b^{r+1}-\b}{\b-1}
\right\vert  \le
D_1\sum_{n\in[1,r]}\left(\frac{\b}{e^{\e }}\right)^n\le
D_2\frac{\b^r}{e^{\e r}}  .  \eqno (3.27)$$
For ${\rm tr}\tr\Psi_\tau$ defined as in (3.19) and $r\ge 1$ we set 
$$\pi^\prime(r)=\sum_{ \lambda(\tau)\le r }
{\rm tr}(\tr\Psi_\tau)    \eqno (3.28)$$
where $\tau$ varies among the prime orbits.

For $\g>1$ we set $r=\g y$; (3.28) implies the equality below; the  inequality follows since 
$\frac{\l(\tau)}{y}> 1$ on the set of the first sum, while $\frac{\l(\tau)}{y}\ge 0$ for all $\tau$; recall that ${\rm tr}\tr\Psi_\tau^k\ge 0$ by Lemma 3.2. The second equality comes from the definition of $\eta$ in (3.26) and the notation (3.20). 
$$\pi^\prime(r)=\pi^\prime(y)+
\sum_{\{ \tau\st y<\l(\tau)\le r \}}{\rm tr}(\tr\Psi_{\tau})\le$$
$$\pi^\prime(y)+\sum_{\{ (k,\tau)\st k\cdot\l(\tau)\le r \}}\frac{\l(\tau)}{y}
{\rm tr}(\tr\Psi_{\tau})^k=
\pi^\prime(y)+\frac{1}{y}\eta(r)  .  $$
If we multiply this inequality by $\frac{r}{\b^r}$ and recall that $r=\g y$, we get 
$$\frac{r\cdot\pi^\prime(r)}{\b^r} \le
\frac{(\g y)\cdot\pi^\prime(y)}{\b^{\g y}}+
\frac{\g\cdot\eta(r)}{\b^r}  .   \eqno (3.29)$$
In Lemma 3.4 below we are going to show that, if $\g^\prime>1$, then 
$$\frac{\pi^\prime(y)}{\b^{\g^\prime y}}\tends 0
\txt{as}y\tends+\infty  .  \eqno (3.30)$$
We show that the last two formulas imply  
$$\limsup_{r\tends+\infty}
\frac{r\cdot\pi^\prime(r)}{\b^r}\le
\g\cdot\frac{\b}{\b-1}  .  \eqno (3.31)$$
For the proof, we fix $\g^\prime\in(1,\g)$; (3.29) implies the first inequality below, while (3.30) and the fact that $\b^{\g-\g^\prime}>1$ imply the second one; the last one comes from (3.27). 
$$\limsup_{r\tends+\infty}\frac{r\cdot\pi^\prime(r)}{\b^r}\le
\g\limsup_{y\tends+\infty}\frac{y}{\b^{(\g-\g^\prime)y}}\cdot
\frac{\pi^\prime(y)}{\b^{\g^\prime y}}  +
\g\limsup_{r\tends+\infty}\frac{\eta(r)}{\b^r}\le$$
$$\g\limsup_{r\tends+\infty}\frac{\eta(r)}{\b^r}\le
\g\cdot\frac{\b}{\b-1}  .  $$
This shows (3.31).

Since $\g>1$ is arbitrary, (3.31) implies that 
$$\limsup_{r\tends+\infty}\frac{\pi^\prime(r)r}{\b^r}\le\frac{\b}{\b-1}  .  $$
Consider now the function $\pi$ of (11) with $\hat V\equiv 1$; in this case it is immediate that the constant $c$ of theorem 2 is $\log\b>0$, where $\b$ is the eigenvalue for the zero potential. By the definition of 
$\pi^\prime$ in (3.28) we get that 
$$\pi(r)=\sum_{e^{c\l(\tau)}\le r}{\rm tr}\Psi_\tau=
\pi^\prime\left(
\frac{\log r}{c}
\right)  .  $$
Since $c=\log\b$, the last two formulas imply (12).

The case $\b<1$  is simpler: (3.26) and (3.28) imply the first inequality below, the second one comes from (3.25) and the last one from the fact that $\b<1$. 
$$\pi^\prime(r)\le\eta(r)\le$$
$$\sum_{n\in[1,r]}\left[
\b^n+D_1\frac{\b^n}{e^{n\e}}
\right]   \le D_2    .  $$
Now (12) follows. 

Now there is only one only missing ingredient, the proof of (3.30); we begin with some preliminaries. 

We note that, by point 5) of Lemma 3.2,
$$\sum_{i=1}^m\a_{\tau,i}^k\ge 0  .  $$
If we assume that $x>0$, this implies the first inequality below; the second inequality comes from the fact that $x\ge\log(1+x)$ and the formula above again. The two equalities that bookend the formula come from the properties of the logarithm. 
$$\prod_{i=1}^m\frac{1}{1-\a_{\tau,i}\cdot x}=
\exp\left[
\sum_{k\ge 1}\frac{1}{k}x^k\left(
\sum_{i=1}^m\a_{\tau,i}^k
\right)
\right] \ge$$
$$\exp\left[
\left(
\sum_{i=1}^m\a_{\tau,i}
\right)x
\right]   \ge  $$
$$\exp\left[
\left(
\sum_{i=1}^m\a_{\tau,i}
\right)\log(1+x)
\right]   =   
\Big( 1+x \Big)^{\sum_{i=1}^m\a_{\tau,i}}   .   \eqno (3.32)$$

\lem{3.4} Formula (3.30) holds for all $\g^\prime>1$.

\proof Since $\g^\prime>1$ and we are supposing $\b>1$, 
$$\frac{1}{\b^{\g^\prime}}<\frac{1}{\b}  .  $$
Together with Proposition 3.3 this implies that $\frac{1}{\b^{\g^\prime}}$ belongs to the disk where  $\zeta(z,0)$ converges. Now we take $y>1$; since $\frac{1}{\b^{\g^\prime}}>0$ and 
${\rm tr}(\Psi_\tau)>0$ by Lemma 3.2, (3.21) implies the first inequality below, while the first equality comes from the properties of the logarithm; recall that $V\equiv 0$. For the second inequality we apply (3.32) to the innermost product. For the third inequality we note that 
$1+\frac{1}{\b^{\g^\prime\cdot y}}$ is smaller than all the factors in the product. The second equality comes from the definition of $\pi^\prime$ in (3.28). The second last inequality follows from the binomial theorem while the last one is obvious. 
$$\zeta(\frac{1}{\b^{\g^\prime}},0)\ge
\exp\sum_{\l(\tau)\le y}\sum_{j=1}^m\sum_{k\ge 1}
\frac{\a^k_{\tau,j}e^{kV(\tau)}\frac{1}{\b^{\g^\prime k\l(\tau)}}}{k}=$$
$$\prod_{\l(\tau)\le y}\prod_{i=1}^m
\frac{1}{1-\frac{\a_{\tau,i}}{\b^{\g^\prime\cdot\l(\tau)}}}  \ge  $$
$$\prod_{\l(\tau)\le y}\left(
1+\frac{1}{\b^{\g^\prime\cdot\l(\tau)}}
\right)^{\sum_{i=1}^m\a_{\tau,i}}  \ge
\prod_{ \l(\tau)\le y }
\left(
1+\frac{1}{\b^{\g^\prime\cdot y}}
\right)^{\sum_{i=1}^m\a_{\tau,i}}=$$
$$\left(
1+\frac{1}{\b^{\g^\prime\cdot y}}
\right)^{\pi^\prime(y) }\ge
1+\frac{\pi^\prime(y)}{\b^{\g^\prime\cdot y}}\ge
\frac{\pi^\prime(y) }{\b^{\g^\prime\cdot y}}  .  $$
Thus, $\frac{\pi^\prime(y)}{\b^{\g^\prime\cdot y}}$ is bounded for all $\g^\prime,y>1$, which implies that 
$$\frac{\pi^\prime(y)}{\b^{\g^\prime\cdot y}}\tends 0\txt{for}y\tends+\infty  $$
for all $\g^\prime>1$.  

\fin

\vskip 2pc

\centerline{\bf \S 4}
\centerline{\bf The weakly mixing case}

\vskip 1pc

In this section we shall work with a potential $\hat V$ such that $\hat V>0$. We begin with the following lemma. 

\lem{4.1} Let $c\in\R$, let $\hat V\in C^{0,\a}(\Sigma,\R)$ be positive and let 
$P(-c\hat V)=\log\b_{-c\hat V}$ be as in the notation after Formula (2.3). Then, there is a unique $c\in\R$ such that 
$$P(-c\hat V)=0  .  \eqno (4.1)$$
Moreover, $c>0$ if $P(0)>0$ and $c<0$ if $P(0)<0$. 

\proof The thesis follows from the intermediate value theorem and the three facts below. 

\noindent 1) The function $\fun{}{c}{P(-c\hat V)}$ is continuous. To show this, we recall that  
$P(-c\hat V)=\log\b_{-c\hat V}$ where $\b_{-c\hat V}$ is the maximal eigenvalue of the Ruelle operator $\L_{-c\hat V}$; thus, it suffices to show that $\fun{}{c}{\b_{-c\hat V}}$ is continuous. Recall that Proposition 1.2 follows from the Perron-Frobenius Theorem 1.1, which says that the eigenvector $Q_{-c\hat V}$ is the fixed point of a contraction; now the theorem of contractions depending on a parameter implies that $Q_{-c\hat V}$ (or, better, its ray) and $\b_{-c\hat V}$ depend continuously on $c$.

\noindent 2) The second fact is that the map $\fun{}{c}{P(-c\hat V)}$ is strictly monotone decreasing. To show this, we set $\e\colon=\min\hat V$ and note that $\e>0$ since 
$\hat V$ is continuous and positive and $\Sigma$ is compact.  Let $c^\prime>c$ and let 
$Q$ be a positive eigenvector of $\L_{-c\hat V}$; then we have the equality below, while the inequality (between matrices, of course) comes from (1.24). 
$$\L_{-c^\prime\hat V}Q\le
e^{-\e(c^\prime-c)}\L_{-c\hat V}Q=
\b_{-c\hat V}\cdot e^{-\e(c^\prime-c)}\cdot Q  .  $$
Iterating, we get the second inequality below, while the first one comes from the fact that 
$\L_{-c^\prime\hat V}$ preserves positive-definiteness.  
$$0\le
\left(
\frac{\L_{-c^\prime\hat V}}{\b_{-c\hat V}}
\right)^n  Q  \le
e^{-\e n(c^\prime-c)}Q\tends 0  .  $$
By point 3) of Proposition 1.2 this implies that 
$$0<\b_{-c^\prime\hat V}<\b_{-c\hat V}   ,   $$
as we wanted. 

\noindent 3) The last fact is  that $P(-c\hat V)\tends\pm\infty$ as $c\tends\mp\infty$; since 
$P(-c\hat V)=\log\b_{-c\hat V}$, this is equivalent to $\b_{-c\hat V}\searrow 0$ as 
$c\nearrow+\infty$ and $\b_{-c\hat V}\nearrow+\infty$ if $c\searrow-\infty$. This follows easily from the formula in point 2) above.

\fin

\noindent{\bf Definitions.}  Having fixed $c$ so that (4.1) holds, we set $V=|c|\hat V\ge 0$; in some of the formulas below we shall have to flip signs according to the sign of $c$. 

We shall count periodic orbits according to a "length" induced by $\hat V$; by Lemma 4.3 below it is equivalent whether we use the length induced by $V$ or by $\hat V$, and we shall use the former.

Let $\fun{g}{\Sigma}{\R}$ be a function; let $g^n$ be as in point 5) of Proposition 1.2 and let $Fix_n$ be as in the definition after (3.10); we set  
$$Z(g)=\sum_{n\ge 1}\frac{1}{n}\sum_{x\in{ Fix}_n}
e^{g^n(x)}\cdot{\rm tr}(\tr\Psi_{x,n})   $$
and 
$$\zeta_{-V}(s)=\exp Z(-sV)=$$
$$\exp\sum_{n\ge 1}\frac{1}{n}\sum_{x\in{ Fix}_n}
e^{-V^n(x)\cdot s}\cdot{\rm tr}(\tr\Psi_{x,n})  .  \eqno (4.2)$$

\vskip 1pc

We must understand for which values of $s$ the series above converges; since by (3.11) 
$\zeta_{-V}(s)=\zeta(1,-sV)$, the radius of convergence of $\zeta(\cdot,-sV)$ must be larger than 1. By Proposition 3.3 this holds if $\b_{-sV}<1$; now recall that $c$ satisfies (4.1); since $V=|c|\hat V$, we get that $\b_{-V}=1$ if $c>0$ and $\b_V=1$ if $c<0$. Let us suppose that $c>0$; since the last lemma implies that the function $\fun{}{s}{\b_{-sV}}$ is monotone decreasing, we get that 
$\zeta_{-V}(s)$ converges if $s>1$. If $c<0$, we see that $\b_{-sV}=\b_{sc\hat V}$; recall that $\b_{-c\hat V}=\b_V=1$ and that $\fun{}{s}{\b_{-sV}}$ is monotone decreasing; all these facts imply that $\zeta_V(s)$ converges if $s>-1$. 

Actually, $\zeta_{-V}(s)$ converges for ${\rm Re}(s)>1$ in one case and in 
${\rm Re}(s)>-1$ in the other one; we skip the easy proof of this fact, which is based on the last paragraph and the fact that ${\rm tr}(\tr\Psi_{x,n})\ge 0$. 

Throughout this section we are going to use the notation of (3.19) and (3.20); now we proceed to find an analogue of (3.21). We could use the fact that 
$\zeta_{-V}(s)=\zeta(1,-sV)$, but we repeat the proof of Section 3. 

The first equality below is the definition (4.2). For the second one we recall from lemma 3.2 that, if $x$ is a periodic orbit of minimal period $n$, then the eigenvalues of 
$\tr\Psi_{x,kn}$ are $\a_{x,i}^k=\a_{\tau,i}^k$ and $V^{nk}(x)=kV^n(x)=kV^n(\tau)$. The third equality comes from (3.19) and the fact that $\tau$, the equivalence class of $x$, has $n$ elements. The last two equalities follow from the properties of the logarithm. 
$$\zeta_{-V}(s)=\exp\sum_{n\ge 1}\frac{1}{n}
\sum_{x\in Fix_n}e^{-V^n(x)\cdot s}{\rm tr}(\tr\Psi_{x,n})=$$
$$\exp\sum_{n\ge 1}\frac{1}{n}\sum_{\{\s^n(x)=x,\quad n\quad least \}}\sum_{k\ge 1}
\sum_{i=1}^{m}
\frac{e^{-kV^n(x)\cdot s}}{k}\cdot\a_{x,i}^k=$$
$$\exp\sum_\tau\sum_{k\ge 1}\sum_{i=1}^{m}\frac{
e^{-kV(\tau)\cdot s}\a_{\tau,i}^k
}{
k
}  =  
\exp\left\{
-\sum_\tau\sum_{i=1}^{m}\log\left[
1-e^{-V(\tau)\cdot s}\cdot\a_{\tau,i}
\right]
\right\}  =  $$
$$\prod_\tau\prod_{i=1}^{m}\frac{
1
}{
1-e^{-V(\tau)\cdot s}\cdot\a_{\tau,i}
}   .    \eqno (4.3)$$

\vskip 1pc

\noindent{\bf Standing hypothesis.} Our standing hypothesis from now on is that, when $c>0$, the function $\zeta_{-V}(s)$ can be extended to a continuous function on 
$\{ {\rm Re}(s)\ge 1 \}\setminus \{ 1 \}$, with a simple pole at $s=1$. When $c<0$ we don't make any hypothesis: it will suffice the fact, which we saw above, that $\zeta_{-V}$ extends to a holomorphic function in ${\rm Re}(s)> -1$. In other words, we are supposing that 
$$\zeta_{-V}(s)=\left\{
\eqalign{
\frac{\tilde\a_+(s)}{(s-1)}&\txt{if}c>0\cr
\tilde\a_-(s)&\txt{if}c<0
}
\right.\eqno (4.4)$$
with $\tilde\a_+$ analytic and non-zero in ${\rm Re}(s)\ge 1$ if $c>0$, and $\tilde\a_-$ analytic and non-zero in ${\rm Re}(s)> -1$ if $c<0$. 

\vskip 1pc

\noindent{\bf Remark.} In the scalar case of [15], the fact that the zeta function has no zeroes on the line ${\rm Re}(s)=1$ has a connection with the dynamics, being equivalent to the fact that the suspension flow by $\hat V$ is topologically mixing. One of the ways to express the topologically mixing property is the following. We define the operator 
$$\fun{V}{C(\Sigma,\C)}{C(\Sigma,\C)}  $$
$$(Vw)(x)=e^{-ia\hat V(x)}w\circ\s(x)  .  $$
We say that the suspension flow (which we leave undefined) is weakly topologically mixing if $Vw=w$ implies that $a=0$ and that $w$ is constant. As we said, we shall not study the connection between this property and the behaviour of the zeta function on the critical line.

\vskip 1pc

Next, we set 
$$N(\tau)=e^{V(\tau)}=e^{|c|\hat V(\tau)}  .  \eqno (4.5)$$
With this notation we can re-write one equality of (4.3) as 
$$\zeta_{-V}(s)=\exp\sum_{k\ge 1}\sum_\tau\sum_{i=1}^{m}
\frac{N(\tau)^{-sk}}{k}\a_{\tau,i}^k  .   $$
As we saw after (4.2), this holds in $Re(s)>1$ if $c>0$ , and in $Re(s)>-1$ if $c<0$. 

Now we take derivatives in two different expressions for $\log\zeta_{-V}(s)$; for the equality on the left we differentiate the formula above, for that on the right we differentiate (4.4) and recall that $\tilde\a_\pm\not=0$; the function $\a_+$ is analytic in an open half-plane containing  $Re(s)\ge 1$ and $\a_-$ is analytic in the open half-plane  $Re(s)> -1$. 
$$-\sum_{k\ge 1}\sum_\tau\sum_{i=1}^{m}
[\log N(\tau)]\cdot N(\tau)^{-sk}\cdot\a_{\tau,i}^k=$$
$$\frac{
\zeta_{-V}^\prime(s)
}{
\zeta_{-V}(s)
}=
\left\{
\eqalign{
-\frac{1}{s-1}+\a_+(s)&\txt{if} c>0\cr
\a_-(s)&\txt{if} c<0  .  
}
\right.        \eqno (4.6)$$
We denote by $\lfloor a\rfloor$ the integer part of $a\ge 0$; recall that we saw at the beginning of this section that $V=|c|\hat V>0$; by (4.5) this implies that 
$\log N(\tau)$ is always positive. As a consequence, for 
$r\ge 1$ we can define 
$n(r,\tau)\ge 0$ by 
$$n(r,\tau)\colon=
\left\lfloor
\frac{
\log r
}{
\log N(\tau)
}
\right\rfloor   .    $$
We define the function $\fun{S}{(1,+\infty)}{\R}$ in the first equality below; the second one comes from the definition of $n(r,\tau)$ above and the fact that $N(\tau)\ge 1$. 
$$S(r)\colon=\sum_{\{ (\tau,n)\st N(\tau)^n\le r\}}\sum_{i=1}^{m}
( \log N(\tau) )\cdot\a_{\tau,i}^n=$$
$$\sum_{\{ \tau\st N(\tau)\le r \}}( \log N(\tau) )\cdot
\sum_{n=1}^{n(r,\tau)}
\left( \sum_{i=1}^{m}  \a_{\tau,i}^n   \right)  .  \eqno (4.7)$$
Note that $S$ is monotone increasing and thus its derivative $\dr S$ is a positive measure; from the first sum above we gather that each orbit $\tau^\prime=(\tau,k)$ contributes to 
$\dr S$ the measure 
$$\sum_{i=1}^{m}
\left(\log N(\tau)\right)
\a_{\tau,i}^k   
\cdot\d_{N(\tau)^k}     $$
where $\d_x$ denotes the Dirac delta centred at $x$. Let now $c>0$; the second equality below is the right hand side of (4.6); the first equality follows from the formula above, the left hand side of (4.6) and the fact that $N(\tau)\ge 1$ by (4.5). The integral is Stieltjes and $Re(s)> 1$.
$$-\int_1^{+\infty}r^{-s}\dr S(r)=
\frac{
\zeta^\prime_{-V}(s)
}{
\zeta_{-V}(s)
}=\frac{-1}{s-1}+\a_+(s)  .  \eqno (4.8)_+$$
When $c<0$ we do the same thing, but using the second expression on the right of (4.6); the equalities hold for $s>-1$. 
$$-\int_{1}^{+\infty}r^{-s}\dr S(r)=
\frac{
\zeta^\prime_{-V}(s)
}{
\zeta_{-V}(s)
}=\a_-(s)  .    \eqno (4.8)_-    $$

We recall the Wiener-Ikehara Tauberian theorem. A proof is in appendix I of [14]; for the whole theory, we refer the reader to [25]. 

\thm{4.2} Assume that the formula below holds, where $A\in\R$, the integral on the left is Stieltjes, $\a$ is analytic in ${\rm Re}(s)>1$ and has a continuous extension to 
${\rm Re}(s)\ge 1$; the function $S$ is monotone increasing.
$$\int_1^{+\infty}r^{-s}\dr S(x)=\frac{A}{s-1}+\a(s)  .  $$
Then,
$$\lim_{s\tends+\infty}\frac{S(s)}{s}=A  .  \eqno (4.9)  $$

\rm

\vskip 1pc

For $N(\tau)$ defined as in (4.5) we set, as in (11), 
$$\pi(r)=\sum_{ N(\tau)\le r }\sum_{i=1}^{m}\a_{\tau,i}    .    \eqno (4.10)$$

\lem{4.3} Let the potentials $V,\hat V>0$ be as at the beginning of this section; let 
$\pi$ be defined as in $(4.10)$ and let 
$$\hat\pi(r)=
\sum_{ e^{\hat V(\tau)}\le r }{\rm tr}\tr\Psi_\tau  .  \eqno (4.11)$$
Then, 
$$\hat\pi(r)=
\pi(r^{|c|})  .    \eqno (4.12)$$
In particular, if
$$\lim_{r\tends+\infty}\frac{\pi(r)\log r}{r}\le A  ,  \eqno (4.13)$$
then 
$$\lim_{r\tends+\infty}\frac{\hat\pi(r)\log r}{r^{|c|}}\le\frac{A}{|c|}  .  \eqno (4.14)$$
Lastly, if we had $\ge$ in (4.13), then we would have $\ge$ also in (4.14). 

\proof We define $N(\tau)$ as in (4.5); the first equality below is comes from the fact that $V=|c|\hat V$ and the second one is (4.5). 
$$e^{\hat V(\tau)}=e^{\frac{1}{|c|}V(\tau)}=N(\tau)^\frac{1}{|c|}   .  $$

The first equality below is the definition of $\hat\pi$ in (4.11); since $|c|>0$ the second one comes from the formula above and the third one is the definition of $\pi$ in $(4.10)$. 
$$\hat\pi(r)=
\sum_{ e^{\hat V(\tau) }\le r }\sum_{i=1}^m\a_{\tau,i}=
\sum_{ e^{ V(\tau) }\le r^{|c|} }\sum_{i=1}^m\a_{\tau,i}=
\pi(r^{|c|})   .  $$

This shows (4.12); now a simple calculation shows that (4.13) implies (4.14). 

\fin

\noindent{\bf Proof of point 2) of Theorem 2.} By Lemma 4.3, formula (13) of Theorem 2 follows if we show (4.13) with $A=1$ when $c>0$ and with $A=0$ when $c<0$; this is what we shall do. 

We distinguish two cases: when $c>0$, $(4.8)_+$ implies that $S$ satisfies the hypotheses of Theorem 4.2 for $A=1$ and thus (4.9) holds, always for $A=1$; a consequence is the asymptotic sign in the formula below, which holds for  $r\tends+\infty$. The equality comes from $(4.7)$. 
$$r\simeq S(r)=
\sum_{N(\tau)\le r}  (\log N(\tau))   \sum_{n=1}^{n(t,\tau)}
\left(   \sum_{i=1}^{m}\a_{\tau,i}^n    \right)      .      \eqno (4.15)   $$
Now we consider $\g>1$ and, for $r\ge 1$, we define $y=y(r)$ by the formula below.
$$1\le y=r^\frac{1}{\g}<r  .  \eqno (4.16)$$
The first equality below comes from the definition of $\pi(r)$ in $(4.10)$; for the first inequality we use point 5) of Lemma 3.2 and the fact that  $\frac{\log N(\tau)}{\log y}\ge 1$ if 
$N(\tau)\ge y$ and $\frac{\log N(\tau)}{\log y}\ge 0$ in all cases; this follows from (4.5) and the fact that $V>0$. The second equality comes from (4.7);  the equality at the end comes by (4.16). 
$$\pi(r)=\pi(y)+
\sum_{ y<N(\tau)\le r }{\rm tr}(\tr\Psi_{\tau})\le$$
$$\pi(y)+\sum_{\{ (\tau,k)\st N(\tau)^k\le r \}}\frac{\log N(\tau)}{\log y}
{\rm tr}(\tr\Psi_{\tau}^k)=$$
$$\pi(y)+\frac{1}{\log y}S(r)=\pi(y)+\frac{\g}{\log r}S(r)  .  $$
If we multiply the last formula by $\frac{\log r}{r}$ and recall (4.16) we get
$$(\log r)\cdot\frac{\pi(r)}{r}\le
\frac{\pi(y)}{y^\g}\cdot\g\cdot\log y+\g\cdot\frac{S(r)}{r}  .  \eqno (4.17)$$
In Lemma 4.2 below we are going to see that, if $\g^\prime>1$, then 
$$\limsup_{y\tends+\infty}\frac{\pi(y)}{y^{\g^\prime}}<+\infty  .  \eqno (4.18)$$
Let now $\g^\prime\in(1,\g)$; Formula (4.17) implies the inequality below, (4.18) and (4.15) imply the equality. 
$$\limsup_{r\tends+\infty}\frac{\pi(r)\log r}{r}\le$$
$$\g\cdot\limsup_{y\tends+\infty}
\frac{\pi(y)}{y^{\g^\prime}}\cdot\frac{\log y}{y^{\g-\g^\prime}}+
\g\cdot\limsup_{r\tends+\infty}\frac{S(r)}{r}=\g  .  $$
Recalling that $\g>1$ is arbitrary, we get the inequality of (4.13) with $A=1$ i.e. Formula (13) of Theorem 2. 

The case $c<0$ is similar. Formula $(4.8)_-$ implies that $S$ satisfies the hypotheses of Theorem 4.2 for $A=0$; in turn, this implies that 
$$\frac{S(r)}{r}\tends 0\txt{as}r\tends+\infty  .  $$
Now the same argument as above implies that 
$$\frac{\pi(r)\log(r)}{r}\tends 0\txt{as}r\tends+\infty  .  $$

Now (4.18) is the only missing ingredient of the proof of Theorem 2. 

\lem{4.4} If $\g^\prime>1$, Formula (4.18) holds.

\proof We show the case $c>0$, since the other one is analogous. We take $\g>1$ and 
$y>1$; recall that we saw after (4.2) that $\zeta_{-V}(s)$ is analytic in $Re(s)> 1$; as a consequence, (4.3) holds for $s=\g$; by point 5) of Lemma 3.2 this yields the first inequality below. As for the equality, it comes from the properties of the logarithm as in (4.3). For the second inequality, we use (3.32). For the third inequality we note that $(1+y^{-\g})$ is the smallest argument; the last inequality follows by the binomial theorem. 
$$\zeta_{-V}(\g)\ge
\exp\sum_{N(\tau)\le y}\sum_{k\ge 1}\sum_{i=1}^m
\frac{e^{-k V(\tau)\cdot\g}\a^k_{\tau,i}}{k}=$$
$$\prod_{N(\tau)\le y}\prod_{i=1}^m\frac{1}{1-e^{-V(\tau)\cdot \g}\a_{\tau,i}}\ge
\prod_{N(\tau)\le y}\prod_{i=1}^m[
1+N(\tau)^{-\g}
]^{\a_{i,\tau}} \ge$$
$$(1+y^{-\g})^{\pi(y)}\ge\frac{\pi(y)}{y^\g}  .  $$
Thus, for $y$ sufficiently large, $\frac{\pi(y)}{y^\g}$ is smaller than $\zeta_{-V}(\g)$. 

\fin

\vskip 2pc
\centerline{\bf References} 


\noindent [1] L. Ambrosio, N. Gigli, G. Savar\'e, Gradient Flows, Birkh\"auser, Basel, 2005.




\noindent [2] U.Bessi, Another point of view on Kusuoka's measure, DCDS, {\bf 41-7}, 3251-3271, 2021.  

\noindent [3] U. Bessi, Families  of Kusuoka-like measures on fractals, preprint. 

\noindent [4] G. Birkhoff, Lattice theory, Third edition, AMS Colloquium Publ., Vol. XXV, AMS, Providence, R. I., 1967. 

\noindent [5] D.-J. Feng, A. K\"aenm\"aki, Equilibrium states for the pressure function for products of matrices, Discrete and Contin. Dyn. Syst., {\bf 30}, 699-708, 2011.












\noindent [6] A. Johansson, A. \"Oberg, M. Pollicott, Ergodic theory of Kusuoka's measures, J. Fractal Geom., {\bf 4}, 185-214, 2017.  

\noindent [7] N. Kajino, Analysis and geometry of the measurable Riemannian structure on the Sierpiski gasket, Contemporary mathematics, {\bf 600}, 91-133, 2013. 

\noindent [8] J. Kelliher, Oseledec's multiplicative ergodic theorem, mimeographed notes. 

\noindent [9] J. Kigami, Analysis on fractals, Cambridge, 2001. 

\noindent [10] J. Kigami, M. L. Lapidus, Weyl's problem for the spectral distribution of Laplacians on P. C. F. self-similar fractals, Commun. Math. Phys., {\bf 158}, 93-125, 1993. 


\noindent [11] S. Kusuoka, Dirichlet forms on fractals and products of random matrices, Publ. Res. Inst. Math. Sci., {\bf 25}, 659-680, 1989. 




\noindent [12] R. Ma\~n\'e, Ergodic theory and differentiable dynamics, Berlin, 1983. 




\noindent [13] I. D. Morris, Ergodic properties of matrix equilibrium state, Ergodic Theory and Dyn. Sys., {38/6}, 2295-2320, 2018.

\noindent [14] W. Parry, M. Pollicott, Zeta functions and the periodic orbit structure of hyperbolic dynamics, Asterisque, {\bf 187-188}, 1990. 

\noindent [15] W. Parry, M. Pollicott, An analogue of the prime number theorem for closed orbits of axiom A flows, Annals of Mathematics, {\bf 118}, 573-591, 1983. 

\noindent [16] M. Piraino, The weak Bernoulli property for matrix equilibrium states Ergodic Theory and Dynamical Systems, {\bf 40-8}, 2219-2238, 2020. 

\noindent [17] V. Rokhlin, On the fundamental ideas of measure theory, Transl. Amer. Math. Soc., {\bf 71}, 1952. 

\noindent [18]  Disintegration into conditional measures: Rokhlin's theorem, notes on M. Viana's website. 

\noindent [19] L. Schwarz, Lectures on disintegration of measures, Tata Institute of Fundamental Research, Bombay, 1975. 

\noindent [20] W. Rudin, Real and complex Analysis, New Delhi, 1983. 


\noindent [21] D. Ruelle, Ergodic theory of differentiable dynamical systems, IHES Publ. Math. {\bf 50}, 27-58, 1979. 

\noindent [22] D. Ruelle, Functional determinants related to dynamical dystems and the thermodynamic formalism, mimeographed notes. 


\noindent [23] M. Viana, Stochastic analysis of deterministic systems, Brazilian Math. Colloquium, 1997, IMPA.


\noindent [24] P. Walters, Ergodic theory, an introduction, Springer, Berlin, 1980. 

\noindent [25] N. Wiener, Tauberian theorems. Annals of Math., {\bf 33}, 1-100, 1932.

\end

We are going to see in (4.4) below that something more holds: if $c>0$, $\zeta_{-V}(s)$ is holomorphic for $Re(s)>1$, if $c<0$ it is holomorphic for $Re(s)>-1$. 

\vskip 1pc

\noindent{\bf Weak mixing.} Now we introduce the weak mixing condition. First of all, we must define the suspension flow. We let
$$\tilde\Sigma_{\hat V}=\{
(x,y)\st x\in\Sigma,,\qquad 0\le y\le\hat V(x)
\}  $$
and we set $(x,\hat V(x))\simeq(\s(x),0)$; we define 
$$\Sigma_{\hat V}=\frac{\tilde\Sigma_{\hat V}}{\simeq}  .  $$
On this space we have a semi-flow for positive times given by 
$$\fun{\s_{\hat V,t}}{(x,y)}{(x,y+t)}  .  $$
Roughly, one should think that $\Sigma$ is a Poincar\'e\ section of the flow; the orbit which starts at $(x,0)$ at time $0$ comes back to the section at time $t=\hat V(x)$. 

We say that $\s_{\hat V}$ is weakly topologically mixing if the following happens. Let us suppose that $w\in C(\Sigma_{\hat V},\C)$ is such that, for some $a>0$,  
$$w\circ\s_{\hat V,t}(x,y)=e^{iat}w(x,y)  
\txt{for all} t\ge 0,\qquad (x,y)\in\Sigma_{\hat V}  .  $$
Then, $a=0$ and $w$ is constant. By Proposition 6.2 of [14], the formula above holds if and only if there is a non-zero $w\in C(\Sigma,\R)$ (but $w\in L^2(G,\kappa_{\hat V})$ suffices) such that 
$$w\circ\s(x)=e^{ia\hat V(x)}\cdot w(x)
\txt{for $\kappa_V$-a.e. $x\in\Sigma$.}  $$

It may happen that, while the map $\fun{\s}{\Sigma}{\Sigma}$ is mixing, the flow 
$\fun{\s_{\hat V}}{\Sigma_{\hat V}\times\R}{\Sigma_{\hat V}}$ is not; we had an example of this in the last section where $\hat V\equiv 1$. 

The arguments before Theorem 6.3 of [14] make an important connection: if $\s_{\hat V}$ is weakly topologically mixing, then the $\zeta$ function of (4.2) has

\noindent 1) a simple pole at $s=1$,

\noindent 2) a continuous extension to the  straight line 
$$\{
s\st {\rm Re}(s)=1
\}  \setminus \{ 1 \}  \eqno (4.3)$$

\noindent 3) and no zeroes on the line above. 

To complete the picture, Proposition 8 of [15] says that, if $\hat V(x)$ depends only on the first two coordinates of $x=(x_0x_1\dots)$, then $\s_{\hat V}$ is {\sl not} topologically weakly mixing if the multiplicative group generated by the numbers
$$\{
e^{\hat V(\tau)}\st \txt{$\tau$ is a periodic orbit}
\}   $$
has rank one. This happens, for instance, in the case of the last section, i.e. when
$\hat V\equiv 1$; in this case, (4.2) implies that $\zeta$ is periodic on the line ${\rm Re}(s)=1$; since there is a pole at $s=1$, $\zeta$ cannot be extended to a continuous function up to the set of (4.3). 

The standing hypothesis of this section is that $\s_{\hat V}$ is weakly topologically mixing, and thus that $\zeta$ can be extended to a continuous function up to the set of (4.3). 

\vskip 1pc